\documentclass{amsart}
\usepackage{amsmath}
\usepackage{amsfonts}
\usepackage{amssymb}
\usepackage{amscd}

\DeclareMathOperator{\Rr}{\mathcal R}

\DeclareMathOperator{\Ob}{\mathcal Ob}
\DeclareMathOperator{\Coinv}{\tilde{S}}

\DeclareMathOperator{\E}{\tilde{\varepsilon}}

\DeclareMathOperator{\A}{\mathcal A}
\DeclareMathOperator{\F}{\mathcal F}

\DeclareMathOperator{\C}{\mathbb C}
\DeclareMathOperator{\Z}{\mathbb Z}
\DeclareMathOperator{\Q}{\mathbb Q}

\DeclareMathOperator{\1}{\underline{1}}
\DeclareMathOperator{\GL}{\bf GL}
\DeclareMathOperator{\U}{\mathcal U}
\DeclareMathOperator{\I}{\bf I}
\DeclareMathOperator{\Gm}{\bf G_m}
\DeclareMathOperator{\Ga}{\bf G_a}
\DeclareMathOperator{\Mn}{\bf M}

\DeclareMathOperator{\Span}{span}

\DeclareMathOperator{\Hom}{Hom}
\DeclareMathOperator{\End}{End}

\DeclareMathOperator{\AutOP}{\Aut^{\otimes,\partial}}
\DeclareMathOperator{\K}{\bf k}
\DeclareMathOperator{\Image}{Im}

\DeclareMathOperator{\Aut}{Aut}
\DeclareMathOperator{\IntHom}{\underline{Hom}}

\DeclareMathOperator{\id}{id}
\DeclareMathOperator{\ev}{ev}
\DeclareMathOperator{\Rep}{\bf Rep}

\DeclareMathOperator{\Cat}{\mathcal{C}}

\DeclareMathOperator{\Char}{char}
\DeclareMathOperator{\ord}{ord}

\DeclareMathOperator{\Seq}{\mathcal{V}}

\theoremstyle{plain}
\newtheorem{theorem}{Theorem}
\newtheorem{lemma}{Lemma}
\newtheorem{proposition}{Proposition}
\newtheorem{corollary}{Corollary}

\theoremstyle{definition}
\newtheorem{definition}{Definition}
\newtheorem{example}{Example}

\theoremstyle{remark}
\newtheorem{remark}{Remark}

\newcommand{\Le}{\leqslant}
\newcommand{\Ge}{\geqslant}

\begin{document}
\title[Tannakian approach to differential algebraic groups]{Tannakian approach to linear differential algebraic groups}

\author{Alexey Ovchinnikov}\thanks{The work was partially supported by NSF Grant CCR-0096842
and by the Russian Foundation for Basic Research, project no. 05-01-00671.}
\address{North Carolina State University\\ Department of Mathematics\\
Raleigh, NC 27695-8205, USA}
\curraddr{University of Illinois at
Chicago, Department of Mathematics, 851 S. Morgan Street, M/C 249, Chicago, IL 60607-7045, USA.}
\email{aiovchin@math.uic.edu} 
\urladdr{http://www.math.uic.edu/\textasciitilde aiovchin/} 

\date\today

\maketitle

\begin{abstract}
Tannaka's Theorem states that a linear algebraic group $G$ is
determined by the category of finite dimensional $G$-modules and
the forgetful functor.
We extend this result to linear differential algebraic groups
by introducing a category corresponding to their representations
and show how this category determines such a group.
\end{abstract}

\section{Introduction}
Given a linear algebraic group $G$, a  rational representation
(or finite dimensional $G$-module) is a finite dimensional vector
space $V$ together with  a morphism $\rho_V:G \to \GL(V)$.
The collection of such objects forms a rigid, abelian, tensor
category and Tannaka's theorem
(\cite[Theorem 1]{Saavedra},\cite[Theorem 2.11]{Deligne},\cite[Theorems 2.5.3 and 2.5.7]{Springer})
states that one can recover the group $G$ as an affine variety
together with the morphisms corresponding to multiplication and
inverse (or equivalently, its coordinate ring and its structure as
a Hopf algebra) from the knowledge of this category $\Rep_G$ and the forgetful
functor from $\Rep_G$ to finite dimensional vector spaces.

In this paper, we consider linear differential algebraic groups (or, shorter, linear $\partial$-$\K$-groups).
These are groups of invertible matrices with entries in a
differential field $\K$ of characteristic $0$ with derivation $\partial$ that are, in
addition, differential varieties, that is, they are defined by the
vanishing of differential polynomials.  A representation of such a
group is a finite dimensional vector space $V$ over $\K$ together
with a differential polynomial morphism from $G$ to $\GL(V)$. If $K$ is a $\partial$-field containing $\K$ then one can talk about $K$-points $G(K)$ of the group $G$. In such a way we obtain a functorial definition of a linear differential algebraic group, which we give in Section~\ref{LinAlgGroupsAsFunctors}. In the preceding sections we view such a group as a Kolchin closed subset of $\U^n$, where $\U$ is
a semi-universal $\partial$-extension of the
ground $\partial$-field $\K$.

The study of these groups and their representations was initiated by
Cassidy in \cite{Cassidy,CassidyRep}.  In this paper, we
introduce differentiation on vector spaces (a ``prolongation'' functor) over $\partial$-fields
such that representations of $G$ correspond to
this construction.  We then show that for a linear
differential algebraic group $G$, the category of objects corresponding
to its representations completely determines $G$ as a differential
variety together with its morphisms for multiplication and inverse,
that is, we show that one can recover its {\em differential}
coordinate ring together with its Hopf algebra and {\em differential}
 structure.

The rest of the paper is organized as follows.  Section~\ref{Basics}
gives formal definitions and properties  of linear differential
algebraic groups.  In Section~\ref{Sequences}, we introduce the 
category $\Seq$ and show how a representation of a linear differential algebraic
group corresponds to an object of $\Seq$. In
Section~\ref{PropertiesSection}, we give various representation
theoretic properties of the objects of $\Seq$ as well as some
consequences (e.g., any representation can be constructed from a
faithful representation using the operations of linear algebra and the prolongation functor).
We then show in Section~\ref{ConstantMatricesSection} how to distinguish linear algebraic groups among
all differential algebraic groups using our
method of differentiating representations.
In Section~\ref{TannakasTheoremSection}, we show how to recover the
group from its associated category.
Although for convenience we do most of the computations in the ordinary case, everything goes
through for the partial differential case. The definition of the corresponding category is given in Section~\ref{ManyParameters}.

We note that the categorical approach to representations of linear
algebraic groups leads to the theory of Tannakian categories.
This theory has found many uses and, in particular, one can develop
the Galois theory of linear differential equations using this
categorical approach.  Recently, a theory of parameterized linear
differential equations has been developed by Cassidy and Singer
\cite{PhyllisMichael} where the Galois groups are linear differential
 algebraic groups.  The category $\Seq$ defined in this paper was
motivated by a desire to give a similar categorical development of
the Galois theory of parameterized differential equations.
This program is now realized in the paper \cite{OvchDeligneMilne}.

\section{Basic definitions}\label{Basics}
\subsection{Differential algebra}
A $\Delta$-ring $R$,
where $\Delta = \{\partial_1,\ldots,\partial_m\}$, is a commutative associative ring with unit $1$ and commuting differentiations $\partial_i : R\to R$ such that
$$
\partial_i(a+b) = \partial_i(a)+\partial_i(b),\quad \partial_i(ab) =
\partial_i(a)b + a\partial_i(b)
$$
for all $a, b \in R$. If $\K$ is a field and
a $\Delta$-ring then $\K$ is called a $\Delta$-field. We restrict ourselves to the case of
$$
\Char\K = 0.
$$ If $\Delta = \{\partial\}$ then
a $\Delta$-field is called a $\partial$-field.
For example, $\Q$ is a $\partial$-field with a unique
possible differentiation (which is the zero one). The field
$\C(t)$ is also a $\partial$-field with $\partial(t) = f,$
and this $f$ can be any rational function in $\C(t).$ For simplicity, we will mostly discuss the case of $\Delta=\{\partial\}$ in this paper and
come back to the general case of $m$ commuting
differentiations in Section~\ref{ManyParameters}.
Let 
$$
\Theta = \left\{\partial^i\:|\: i\in \Z_{\Ge 0}\right\}.
$$
Since $\partial$ acts on a $\partial$-ring $R$,
there is a natural action of $\Theta$ on $R$.

A non-commutative ring $R[\partial]$ of linear differential operators is generated as a left $R$-module by the monoid $\Theta$. A typical element
of $R[\partial]$ is a polynomial 
$$
D = \sum_{i=1}^na_i\partial^{i},\ a_i \in R.
$$
The right $R$-module structure follows from the
formula 
$$
\partial\cdot a = a\cdot \partial + \partial(a)
$$
for all $a \in R$.
We denote the set of operators in $R[\partial]$ of order less than
or equal to $p$ by $R[\partial]_{\Le p}.$

Let $R$ be a $\partial$-ring. If $B$ is an $R$-algebra, then $B$ is a $\partial$-$R$-algebra
if the action of $\partial$ on $B$ extends the
action of $\partial$ on $R$. If $R_1$ and $R_2$
are $\partial$-rings then a ring homomorphism
$\varphi: R_1 \to R_2$ is called a $\partial$-homomorphism if it commutes with $\partial$, that is,
$$
\varphi\circ\partial = \partial\circ\varphi.
$$ We denote these homomorphisms simply by $\Hom(R_1,R_2)$.
If $A_1$ and $A_2$ are $\partial$-$\K$-algebras
then a $\partial$-$\K$-homomorhism simply means a 
$\K[\partial]$-homomorphism.
Let $Y = \{y_1,\ldots,y_n\}$ be a set of variables. We differentiate them:
$$
\Theta Y := \left\{\partial^iy_j
\:\big|\: i \in \mathbb{Z}_{\Ge 0},\ 1\Le j\Le n\right\}.
$$
The ring of differential polynomials $R\{Y\}$ in
differential indeterminates $Y$ 
over a $\partial$-ring $R$ is
the ring of commutative polynomials $R[\Theta Y]$
in infinitely many algebraically independent variables $\Theta Y$ with
the differentiation $\partial$, which naturally
extends $\partial$-action on $R$ as follows:
$$
\partial\left(\partial^i y_j\right) := \partial^{i+1}y_j
$$
for all $i \in \Z_{\Ge 0}$ and $1 \Le j \Le n$.
A $\partial$-$\K$-algebra $A$ is called finitely
$\partial$-generated over $\K$ if there exists
a finite subset $X = \{x_1,\ldots,x_n\} \subset A$
such that $A$ is a $\K$-algebra generated by
$\Theta X$.  

An ideal $I$ in a $\partial$-ring $R$ is called differential if it is stable under the action of
$\partial$, that is,
$$
\partial(a) \in I
$$
for all $a \in I$. If $F \subset R$ then $[F]$ denotes the differential ideal generated by $F$.
If a differential ideal is radical, it is called
radical differential ideal. The radical differential
ideal generated by $F$ is denoted by $\{F\}$. If a
differential ideal is prime, it is called a prime
differential ideal.

\subsection{Linear differential algebraic groups}
We shall recall some definitions and results from differential
algebra (see for more detailed information \cite{Cassidy,Kol})
leading up to the ``classical definition'' of a linear differential
algebraic group and its representative functions. Later in the paper
we will give a Hopf-theoretic treatment and provide an equivalent
definition in terms of representable functors.

Let $\K \subset \U$ be a semi-universal differential field over $\K,$
that is, a differential field such that if $K$ is a differential
field extension of $\K,$ finitely generated in the differential
sense, then there exists a $\K$-isomorphism of $K$ into $\U.$ We
will assume that all differential fields we consider are subfields
of $\U$ (the Hopf-theoretic treatment will {\it not} use semi-universal
differential extensions).

\begin{definition} For a differential field extension $K\supset \K$ a {\it Kolchin closed} subset $W(K)$ of $K^n$ over $\K$ is the set of common zeroes
of a system of differential algebraic equations with coefficients in $\K,$ that is, for $f_1,\ldots,f_k \in \K\{Y\}$ we define
$$
W(K) = \left\{ a \in K^n\:|\: f_1(a)=\ldots=f_k(a) = 0\right\}.
$$
\end{definition}

There is a bijective correspondence between Kolchin closed subsets
$W$ of $\U^n$ defined over $\K$ and radical differential ideals
$\I(W) \subset \K\{y_1,\ldots,y_n\}$ 
generated by the differential polynomials $f_1,\ldots,f_k$ that define $W$. 
In fact, the $\partial$-ring $\K\{Y\}$ is
Ritt-Noetherian, meaning that every radical
differential ideal is the radical of a finitely
generated differential ideal, by the Ritt-Raudenbush basis theorem.
Given a Kolchin closed subset $W$ of
$\U^n$ defined over $\K$ we let the {\it coordinate ring} $\K\{W\}$
be:
$$
\K\{W\} = k\{y_1,\ldots,y_n\}\big/\I(W).
$$
A {\it differential polynomial} map $\varphi : W_1\to W_2$ between Kolchin closed subsets of $\U^n,$ defined over $\K,$ is given in coordinates by differential polynomials from $\K\{y_1,\ldots,y_n\}$. To give $\varphi : W_1 \to W_2$
is equivalent to defining $\varphi^* : \K\{W_2\} \to \K\{W_1\}.$

\begin{definition}\cite[Chapter II, Section 1, page 905]{Cassidy} A {\it linear differential
algebraic group} (or linear $\partial$-$\K$-group) is a Kolchin closed subgroup $G$ of $\GL_n(\U),$
that is, an intersection
of a Kolchin closed subset of $\U^{n^2}$ with $\GL_n(\U)$, which is closed under
the group operations.
\end{definition}

Note that we identify $\GL_n(\U)$ with a Zarisky closed
subset of $\U^{n^2+1}$ given by
$$\left\{(A,a)\:\big|\: (\det(A))\cdot a-1=0\right\}.$$
If $X$ is an invertible $n\times n$ matrix, we
can identify it with the pair $(X,1/\det(X))$. Hence, we may represent the coordinate ring of $\GL_n(\U)$ as 
$$
\K\{X,1/\det(X)\}.$$ 
Denote $\GL_1$ simply by $\Gm$. Its coordinate ring
is $\K\{y,1/y\}$, where $y$ is a differential indeterminate.

\begin{definition}\cite{CassidyRep} A differential polynomial
group homomorphism  $\rho : G \to \GL(V)$ is called a
{\it differential representation} of a linear differential algebraic group $G,$ where $V$ is a
finite dimensional vector space over $\K$. Such space is
called a {\it differential $G$-module}.
\end{definition}

A Hopf algebra $A$ is a commutative associative algebra together with
comultiplication $\Delta: A\to A\otimes A$, coinverse $S: A\to A$, and counit $\varepsilon : A \to \K$ satisfying certain axioms \cite[2.1.2]{Springer}.
In \cite[page 23]{Water} representations of an algebraic group $G$
are viewed as comodules over the Hopf algebra $A$ of regular
functions on $G$. In Section~\ref{Sequences}
we will develop a similar technique to look at differential representations
as {\it differential comodules} (Definition~\ref{DiffComodules}) over {\it differential}
Hopf algebras (Definition~\ref{DiffHopfAlgebra}).

Following \cite[page 5]{Water} one interprets algebraic groups as representable
functors from the category of $\K$-algebras to groups, that is, there
must be a $\K$-algebra $A$ such that for any $\K$-algebra $B$ the $B$-points
of $G,$ denoted by $G(B),$ are just $\Hom(A,B).$ Then \cite[Theorem, page 6]{Water} says that natural maps from one group $G$ to another one $G'$, viewed
as functors, correspond to the algebra homomorphisms $A' \to A$.
We will develop this in differential setting in Section~\ref{LinAlgGroupsAsFunctors} not assuming
that all $\partial$-fields we consider are subfields of $\U$.

\subsection{Representative functions}
Let $W \subset \U^n$ be a Kolchin closed subset defined over $\K.$
We define $\K\langle W\rangle$ to be the complete ring of quotients
of $\K\{W\},$ that is $\K\{W\}$ is localized with respect to
the set of nonzero divisors. We note that $\K\langle W\rangle$ is clearly a $\partial$-ring.

A linear differential algebraic group $G(\U)$ acts on the rational
functions $\U\langle G\rangle = \U\otimes\K\langle G\rangle$ by
right translations
(see \cite[page 901]{Cassidy}, \cite[page 227]{CassidyRep}). We will define this independently of $\U$ later on in formula~\eqref{RegularAction} of Section~\ref{RecoveringSection}. According to \cite[page 227]{CassidyRep}
the functions whose orbit generates a finite dimensional vector space are called {\it representative functions}. Note that
the representative functions form a $\partial$-$\K$-algebra, which we denote by $R(G)$.
By  (\cite[page 230, Theorem]{CassidyRep}), $R(G) = \U\{G\}$.

We show how these functions (and, hence, the algebra) can be connected
with finite dimensional differential representations of $G$. Let $V$ be a vector space
over $\K$ of dimension $n$ and let $\rho : G\to \GL(V)$ be a
representation of $G$. The image group $H$ has
coordinate ring $\K\{H\}$. The $\partial$-$\K$-algebra $\K\{H\}$ is a quotient of $\K\{\GL_n\} = \K\{X,1/\det(X)\}$. Let $Z$ be the image of $X$
in $\K\{H\}$ with respect to this canonical homormorphism. Then, $\rho$ induces a $\partial$-$\K$-homomorphism $\rho^* : \K\{H\}\to \K\{G\}$,
mapping the entries of $Z$ onto elements $\varphi_{ij}$ called {\it coordinate functions}
of $\rho$.

\begin{proposition}\label{RepresentableFunctions} Representative functions are the
same as the coordinate functions of finite dimensional representations of $G$.
\end{proposition}
\begin{proof} Let $\rho : G \to \GL(V)$.  Let also $\varphi_{ij}$ be a coordinate function of $\GL(V)$. We have
\begin{equation*}
\begin{CD}
G @>{\rho}>> \GL(V) @>{\varphi_{ij}}>> \U
\end{CD}
\end{equation*}
According to \cite[Corollary 1, page 231]{CassidyRep} we have
$\varphi_{ij}\circ \rho \in \U\{G\}$. By \cite[Theorem, page 230]{CassidyRep}
we have
$R(G) = \U\{G\}$. We show that any $f \in R(G)$ is of the form
$\varphi_{ij}\circ \rho_f$ for a finite dimensional
representation $\rho_f$ of $G$.

Take any $f \in R(G)$ and consider its $G$-orbit
$Gf =: V$,
which is finite dimensional, under the
action $\rho : G \to \GL(\U\{G\})$ with
$$\rho(g)(f)(x) = f(xg).$$
Let $V$ be spanned by $\{f=f_1,\ldots,f_n\}$ over $\U$.
So, $$\rho(g)(f_1) = \sum_{i = 1}^n c_i(g)f_i.$$
Evaluating the last equality at the point
$e \in G$ for all $g \in G$ we get
$$
f(g) = (\rho(g)(f))(e) =  \sum_{i=1}^n c_i(g)\cdot f_i(e).
$$
It remains to change the basis which correspond to the conjugation
of the representation matrices.
Assume that $f_1(e) \ne 0$. A conjugation matrix is
$$
C = \begin{pmatrix}
f_1(e)&f_2(e)& f_3(e) &\ldots&f_n(e)\\
0& 1 & 0&\cdots&0\\
0& 0 & 1&\cdots &0\\
\vdots &\vdots&\vdots &\ddots&\vdots\\
0 & 0 & 0 &\cdots& 1
\end{pmatrix}.
$$
If $B(g)$ is a representation matrix then it goes to $CB(g)C^{-1}$
under such a change of coordinates.
Thus, we have obtained the function $f(g)$ as a coordinate function
of some finite dimensional differential representation of the
linear differential algebraic group $G$.
\end{proof}

In the following we will eliminate semi-universal
$\partial$-extensions of $\K$ by using a functorial approach.

\section{The category $\Seq$ and representations}\label{Sequences}
We start with introducing a {\it differentiation  functor} on finite dimensional
vector spaces and investigate its essential properties.
\subsection{Definition}

\begin{definition}\label{SeqCategory} The category $\Seq$ over a
$\partial$-field $\K$ is the category of
finite dimensional vector spaces over $\K$:
\begin{enumerate}
\item objects are finite dimensional $\K$-vector spaces, 
\item morphisms are $\K$-linear maps;
\end{enumerate}
with tensor product
$\otimes,$ direct sum $\oplus,$ dual $*,$ and additional operations:
$$
\partial^p : V \mapsto V^{(p)} := \K[\partial]_{\Le p}\otimes V,
$$
which we call differentiation (or prolongation) functors. If $\varphi \in \Hom(V,W)$ then we define
$$
\partial^p(\varphi) : V^{(p)}\to W^{(p)},\ \varphi(\partial^q\otimes v) = \partial^q\otimes\varphi(v),\ 0\Le q\Le p.
$$
\end{definition}

\begin{remark} Note that $\K[\partial]$ is a non-commutative $\K$-algebra, so
in this tensor product we think of $\K[\partial]$ as a right $\K$-space and
$V$ as a left $\K$-space.
\end{remark}

\begin{remark}
If $\{v_1,\ldots,v_n\}$ is a basis of $V$ then
$
\{v_1,\ldots,v_n,\ldots,\partial^p\otimes v_1,\ldots,\partial^p\otimes v_n\}
$
is a basis of $V^{(p)}.$
\end{remark}

We denote $\partial \otimes v$ simply by $\partial v.$

\subsection{Linear differential algebraic groups}\label{LinAlgGroupsAsFunctors}
In the following Section~\ref{Representations} we will introduce another definition of a differential
representation of a linear differential algebraic group
(see Definition~\ref{DiffComodules}). For this purpose we will view linear
differential algebraic groups as {\it representable functors} from the
category of $\partial$-$\K$-algebras to groups (see \eqref{GroupsAsFunctors}). So, throughout this
section we are defining a linear differential algebraic group not as a subgroup of $\GL_n(\U)$
but functorially.

We take a Hopf-theoretic approach to the study of linear differential
algebraic groups as in \cite{CassidyRep}. Let $A$ be a
(finitely generated) $\partial$-$\K$-algebra.
Following \cite[page 226]{CassidyRep} one defines
the set $$G(\U) = \Hom(A, \U),$$
where $\U$ is the semi-universal differential field as before, to get $G(\U)$ back from $A$ as
a Kolchin closed subset of $\U^n$.
Assume that $A$ is supplied with the following operations:
\begin{itemize}
\item differential algebra homomorphism $m : A\otimes A \to A$ is the
multiplication map on $A,$
\item differential algebra homomorphism $\Delta : A \to A\otimes A$,
which is a comultiplication, 
\item differential algebra homomorphism $\varepsilon : A  \to k$,
which is a counit,
\item differential algebra homomorphism $S : A \to A$,
which is a coinverse.
\end{itemize}
We also assume that these maps satisfy commutative diagrams (see \cite[page 225]{CassidyRep}):
\begin{equation}\label{HopfAxioms}
\begin{CD}
A @>{\Delta}>> A\otimes A\\
@V{\Delta}VV @V{\id_A\otimes\Delta}VV \\
A\otimes A @>{\Delta\otimes \id_A}>> A\otimes A\otimes A
\end{CD}\qquad
\begin{CD}
A @>{\Delta}>> A\otimes A\\
@V{\id_A}VV @V{\id_A\otimes\varepsilon}VV \\
A @>{\sim}>> A\otimes k
\end{CD}\quad\quad
\begin{CD}
A @>{\Delta}>> A\otimes A\\
@V{\varepsilon}VV @V{m\circ(S\otimes \id_A)}VV \\
k @>{\hookrightarrow}>> A
\end{CD}
\end{equation}

\begin{definition}\label{DiffHopfAlgebra} Such a commutative associative $\partial$-$\K$-algebra $A$
with the unity and operations $m,$ $\Delta,$ $S,$
and $\varepsilon$ satisfying axioms~\eqref{HopfAxioms} is called a {\it differential Hopf algebra} (or {\it Hopf $\partial$-$\K$-algebra}).
\end{definition}

\begin{remark} Introduced in this way $G(\U)$ with operations
corresponding to $\Delta,$ $S,$ and $\varepsilon$ is not only a Kolchin-closed
subset of some $\U^n$ but a group at the same time. Indeed, the group multiplication is defined
in the usual way:
$$
(\varphi_1\cdot\varphi_2)(a) := (\varphi_1\otimes\varphi_2)(\Delta(a))
$$
for all $\varphi_1,\varphi_2 \in \Hom(A,\U) = G(\U)$ and $a\in A$; the group inverse is given by a similar
formula. 
Moreover,
all these operations are continuous in the Kolchin topology (see \cite{Cassidy} for proofs). Finally, since $A$ is $\partial$-finitely generated, the differential algebraic group $G(\U)$ is linear \cite[Proposition 12, page 914]{Cassidy}.
\end{remark}

Recall that a linear $\partial$-$\K$-group $G$ 
is defined  by a system
of differential polynomial equations $F = 0$ 
with coefficients in $\K$ or by the radical
differential ideal $I$ of $\K\{y_1,\ldots,y_n\}$ generated by $F$
(see \cite[page 895]{Cassidy}).
If we represent $0\to I \to \K\{y_1,\ldots,y_n\} \to A \to 0$
then we can consider the group in a $\partial$-$\K$-algebra $B,$ that
is,
$$
G(B) = \Hom_{\K[\partial]}(A,B) =: \Hom(A,B).
$$
Note that here $B$ is not necessarily a subring of
$\U$.
For convenience we gather all such $G(B)$ and form the
following {\it representable} functor
\begin{align}\label{GroupsAsFunctors}
G : \{\text{$\partial$-$\K$-algebras}\} \to \{\text{Groups}\}, \quad B \mapsto \Hom_{\K[\partial]}(A,B),
\end{align}
which we call the {\it affine differential algebraic group} defined by $A$.

Among all differential algebraic groups we distinguish the differential
general linear group. Consider a finite dimensional vector space $V$
over $\K$ of dimension $m$. We define $\GL\left(V\right)$ by the functor
$$
\GL\left(V\right) : B \mapsto \Hom_{\K[\partial]}(\K\{X_{11},\ldots,X_{mm},1/\det\},B).
$$
One may consider $\GL\left(V\right)(B)$ as the set of $m\times m$ invertible
matrices with coefficients in the $\partial$-$\K$-algebra $B$.
\begin{example}
Recall, that the coordinate ring of $\Gm$ is
$\K\{y,1/y\}$. Its $\partial$-$\K$-Hopf algebra
structure is given by 
\begin{align*}
\Delta(y) &= y\otimes y,\\
S(y) &= 1/y.
\end{align*}
These maps are $\partial$-homomorphisms. Therefore,
\begin{align*}
\Delta(\partial y) = \partial(\Delta(y))=\partial y\otimes y + y\otimes\partial y,\\
S(\partial y) = \partial(S(y)) = \partial(1/y) = -\partial y/y^2
\end{align*}
By differentiating these expressions further
one gets the action of $\Delta$ and $S$ on higher derivatives of $y$.
\end{example}
\begin{example}The differential additive group $\Ga$ is represented by the $\partial$-$\K$-Hopf algebra $\K\{y\}$ 
with 
\begin{align*}
\Delta(\partial^py) &= \partial^py\otimes 1 + 1\otimes\partial^py,\\
S(\partial^py) &= (-1)^{p+1}y
\end{align*}
for all $p \in\Z_{\Ge 0}$.
\end{example}
The linear differential algebraic groups we defined earlier correspond
to subgroups of $\GL(V).$ Unlike the situation for algebraic groups,
there are affine differential algebraic groups that are not isomorphic
to linear differential algebraic groups (see \cite[page 911]{Cassidy}).

{\it Morphisms of differential algebraic groups} then in our sense are {\it morphisms of
their representable functors}.
We need the following result (a corollary of Yoneda's Lemma) from the theory of categories to
see that this corresponds to the morphisms of algebras defining the groups.

\begin{lemma}\cite[Corollary 2, page 44]{Pareigis},\cite[30.7, Corollary, page 224]{Herrlich}\label{Yoneda}
Let $\Cat$ be a category such that $\Hom(A,B)$ is a set for all
objects $A, B$ of $\Cat$.
Let $E$ and $F$ be functors from the category $\Cat$ to the category of sets represented
by some objects $A$ and $B$, that is, $E = \Hom\left(A,\_\:\right)$ and
$F = \Hom\left(B,\_\:\right)$. Then morphisms of functors $E$ and $F$
correspond to homomorphisms of $B$ and $A$.
\end{lemma}

It remains to note that the category of $\partial$-$\K$-algebras
satisfies the assumptions of Lemma~\ref{Yoneda}.

\subsection{Representations}\label{Representations}
We will take a careful look at differential representations of
a linear differential algebraic group $G$. These are differential
algebraic group homomorphisms
$$
\Phi : G \to \GL\left(V\right) =: \Aut\left(V\right)
$$
for some finite dimensional
$\K$-vector space $V$. 
\begin{remark}Here, $G$ and $\GL(V)$
are considered as functors whose points vary with
$\partial$-$\K$-algebras as it was explained in
Section~\ref{LinAlgGroupsAsFunctors}. As mentioned
earlier, $G$ is determined by its $\U$-points and
we can {\it identify} the two concepts: $G(\U)$ and $G$.
If $\K$ is a differentially closed field then
$G(\K)$ determines $G$ and one does not need to
look at $G(\U)$.
\end{remark}
By Lemma~\ref{Yoneda} the morphism $\Phi$  corresponds to
the homomorphism of the $\partial$-$\K$-algebras.
In \cite{CassidyRep} a representation of the group $G$
defined over $\K$ is
a rational differential algebraic group homomorphism
$$
G(\U) \to \GL\left(V\right)(\U).
$$
But such a morphism is a differential polynomial map
by \cite[Corollary 1, page 231]{CassidyRep} and so corresponds to a homomorphism
of the associated $\partial$-$\K$-algebras. Hence, we
can freely use our language of functors.

\begin{definition}\label{DiffModules} For a linear differential algebraic
group $G$ an object  $V \in \Ob(\Seq)$ together with a natural map of
group functors $r_V : G \to \Aut(V)$ which is a group homomorphism is called
a {\it differential $G$-module}. The map $r_V$ is called a $\partial$-$\K$-representation of $G$ in $V$. 
\end{definition}

\subsection{Differential comodules}
We are going to restate this in the language of comodules which we introduce
now. For this we define a differential analogue of an algebraic comodule. Let $A$ be a
differential Hopf algebra.

\begin{definition}\label{DiffComodules} A finite dimensional vector space $V$ over $\K$
is called an {\it $A$-differential comodule}
if there is a given $\K$-linear morphism
$$
\rho : V \to V\otimes A
$$
satisfying the axioms:
\begin{equation*}
\begin{CD}
V@>{\rho}>> V\otimes A\\
@VV{\rho}V @VV{\id_V\otimes\Delta}V \\
V\otimes A@>{\rho\otimes\id_A}>> V\otimes A\otimes A
\end{CD}\qquad\qquad
\begin{CD}
V @>{\rho}>> V\otimes A\\
@VV{\id_V}V @VV{\id_V\otimes\varepsilon}V \\
V @>{\sim}>> V\otimes\K
\end{CD}
\end{equation*}
together with the prolongation of $\rho$ on $V^{(i)}$ commuting with
$\partial.$
\end{definition}
The definition and correctness of the prolongation are given in
Theorem~\ref{EquivalentDefsOfRepresentations} and Lemma~\ref{IthRestriction}, which follows the
theorem.
We will show that $A$-differential comodules are in one-to-one
correspondence with differential $G$-modules, where $G$ is the functor
represented by $A$.

\subsection{Equivalent definitions of differential representations}

\begin{theorem}\label{EquivalentDefsOfRepresentations}
Let $A$ be a $\partial$-$\K$-Hopf algebra and $G$
be the linear $\partial$-$\K$-group ($\partial$-$\K$-group functor) represented by $A$. Let $V$
be an object in $\Seq$. Then, there is a bijective correspondence between the set of
$\partial$-representations from $G$ into $\GL(V)$ and the set $\partial$-$A$-comodule structures
$$\rho : V \to V\otimes A$$
on $V$.
If $\{v_1,\ldots,v_n\}$ is a basis of $V$ then in coordinates we have
$$
\rho (v_j) = \sum_{i=1}^nv_i\otimes a_{ij}, \quad \Delta(a_{ij}) = \sum_{r=1}^na_{ir}\otimes a_{rj}.
$$
Moreover, $$\rho(\partial ^pv_j)= \sum_{i=1}^n\sum_{q=0}^p\binom{p}{q}\partial ^q v_i\otimes(\partial^{p-q} a_{ij})$$ gives a prolongation
of $\rho$ on $V^{(p)}$ for all $p\in\Z_{\Ge 1}$.
\end{theorem}
\begin{proof}
The representation $\Phi$ is a morphism of group functors $G \to \Aut(V).$
According to \cite[Theorem, Section 3.2]{Water} such
a representation $\Phi$ of the affine group scheme $G$ corresponds
to the $\K$-linear map $\rho$ with the following commuting diagrams:
\begin{equation*}
\begin{CD}
V @>{\rho}>> V\otimes A\\
@VV{\rho}V @VV{\id_V\otimes\Delta}V \\
V\otimes A @>{\rho\otimes\id_A}>> V\otimes A\otimes A
\end{CD}\qquad\qquad
\begin{CD}
V @>{\rho}>> V\otimes A\\
@VV{\id_V}V @VV{\id_V\otimes \varepsilon}V \\
V @>{\sim}>>  V\otimes\K
\end{CD}
\end{equation*}
More precisely, the map $\rho$ comes from the restriction of the $A$-linear map
$$\Phi(\id_A) : V\otimes A \to V\otimes A$$ to $V\otimes\K \cong V$.
By \cite[Corollary, Section 3.2]{Water} we have
$$
\rho (v_j) = \sum_{i=1}^n v_i\otimes a_{ij}, \quad \Delta(a_{ij})
= \sum_{r=1}^na_{ir}\otimes a_{rj}.
$$
Let us demonstrate the last differential identity. We have
$$\Phi(\id_A) \in \Hom_A(V\otimes A, V\otimes A)$$
and $\Phi(\id_A)$ can be extended to a map $V^{(p)}\otimes A \to V^{(p)}\otimes A$ commuting with the $\partial $-structure.
The first step gives us the following:
\begin{equation*}
\begin{CD}
v_j\otimes 1@>{{\partial }_{V\otimes A}}>> \partial v_j\otimes 1\\
@VV{\Phi(\id_A) = \rho}V @VV{\Phi(\id_A) = \rho}V \\
\sum\limits_{i=1}^nv_i\otimes a_{ij} @>{{\partial }_{V\otimes A}}>>
\sum\limits_{i=1}^n((\partial v_i)\otimes a_{ij}+ v_i\otimes\partial a_{ij})
\end{CD}\qquad\qquad
\end{equation*}
The formula for higher order derivatives can be obtained by induction.
This, indeed, makes $V^{(p)}$ an $A$-comodule (see Lemma~\ref{IthRestriction}).

On the other hand, having such a $\rho : V \to V\otimes A$ one extends
it by the $A$-linearity to $\rho_A : V\otimes A \to V\otimes A$ and  then to $V^{(p)}\otimes A \to V^{(p)}\otimes A$ commuting
with the $\partial $-structure (see Lemma~\ref{IthRestriction} for
correctness). Consider a $\partial$-$\K$-algebra $B$. For
any $g \in \Hom_{\K[\partial]}(A, B)$ in $G(B)$ we have the following commutative diagram:
\begin{equation*}
\begin{CD}
V\otimes A @>{\Phi(\id_A)}>> V\otimes A\\
@VV{\id_V\otimes g}V @VV{\id_V\otimes g}V \\
V\otimes B @>{\Phi(g)}>> V\otimes B
\end{CD}\qquad\qquad
\end{equation*}
meaning that $\Phi(g)$ is determined by $(\id_V\otimes g)\circ\rho$ using
$B$-linearity of $\Phi(g),$ where $\rho$ is the restriction of $\rho_A := \Phi(\id_A)$
to $V.$
Since $\Phi(\id_A)$ is constructed so as to preserve the $\partial $-structure, the map $\Phi(g)$ does the same thing, since
$g$ is a $\partial$-$\K$-algebra homomorphism $A\to B$.
Indeed, take $v = \sum c_i\partial ^iv_i$.
Then
\begin{align*}
\Phi(g)(v) &= (\id_V\otimes g)\circ\Phi(\id_A)\left(\sum c_i\partial ^iv_i\right)=
\sum c_i\cdot(\id_V\otimes g)\partial ^i(\Phi(\id_A)(v_i))=\\
&=\sum c_i\cdot(\id_V\otimes g)\partial ^i
\left(\sum_{j=1}^n v_j\otimes b_{ji}\right)=\\
&=\sum c_i\cdot(\id_V\otimes g)\left(\sum_{j=1}^n
\sum_{r=0}^i \binom{i}{r}\partial ^rv_j\otimes\partial^{i-r}b_{ji}\right)=\\
&=\sum c_i\left(\sum_{j=1}^n\sum_{r=0}^i \binom{i}{r}\partial^rv_j\otimes \partial^{i-r}g(b_{ji})\right)=\\
&=\sum c_i\partial ^i\left(\sum_{j=1}^n v_j\otimes g(b_{ji})\right)= \sum c_i\partial ^i(\Phi(g)(v_i)).
\end{align*}
Here we denoted $\Phi(\id_A)(v_i) = \sum_j v_j\otimes b_{ji}$ for some
elements $b_{ij} \in A.$ From this it follows that
$$
\Phi(g) \in \Hom_{\K[\partial]}\left(\K\left\{X_{11},\ldots,X_{nn}, 1/\det\right\},B\right).
$$
This finally establishes a bijection between differential representations
and differential comodules.
\end{proof}

\subsection{Prolongation of representations}
Let $G$ be a linear differential algebraic group
represented by $A$ and
$\rho : G \to \GL\left(V\right)$ be its differential representation.
\begin{lemma}\label{IthRestriction} The action of $G$ on each ${V^{(i)}}$ is
algebraic, that is, ${V^{(i)}}$ is an $A$-comodule for all $i\Ge 0$.
\end{lemma}
\begin{proof} We have a differential algebraic action
on $V$. Let $\{v_1,\ldots,e_v\}$ be a $\K$-basis of $V.$
For $\rho : V \to V\otimes A$ from Theorem~\ref{EquivalentDefsOfRepresentations} we have
\begin{align*}
\rho(v_j) &= \sum_{i=1}^n v_i\otimes a_{i,j},\\
\Delta(a_{i,j}) &= \sum_{r=1}^n a_{i,r}\otimes a_{r,j},\\
\Delta(\partial^p a_{i,j})=\partial^p(\Delta(a_{i,j})) &=
\sum_{r=1}^n\sum_{q=0}^p\binom{p}{q} (\partial^q a_{i,r})\otimes (\partial^{p-q} a_{r,j}),\\
\rho(\partial ^pv_j) = \partial ^p\rho(v_j) &= \sum_{i=1}^n\sum_{q=0}^p \binom{p}{q}(\partial ^qv_i)\otimes(\partial^{p-q}a_{i,j}).\\
\end{align*}
One can show by induction that if $B = (a_{i,j})_{i,j=1}^n$ is
the ``representation'' matrix for $G \to \GL\left(V\right)$ then for a fixed number $i \Ge 0$ we have
$$
B_i := \begin{pmatrix}
B &0&0&\ldots &0 \\
\binom{i}{1}\cdot B_t &B& 0& \ldots & 0 \\
\binom{i}{2}\cdot B_{tt}&\binom{i-1}{1}\cdot B_t&B&\ldots &0\\
\cdots&\cdots&\cdots&\cdots&\cdots\\
B_{t^i}&\ldots&\ldots&\ldots&B
\end{pmatrix} = (c_{r,s})_{r,s=1}^{i\cdot n}
$$ is the one for $G \to \GL\left({V^{(i)}}\right)$, where $B_{t^k}$ means the
matrix $\left(\partial^ka_{i,j}\right)_{i,j=1}^n.$

It remains to show we do have a group action. Let us do this.
The ordered basis of ${V^{(i)}}$ is $$\{\partial ^iv_1,\ldots,\partial ^iv_n,\partial ^{i-1}v_1,
\ldots,\partial ^{i-1}v_n,\ldots,v_1,\ldots,v_n\}.$$ By induction,
using the facts that $B_i$ has:
\begin{enumerate}
\item $B$ on the main diagonal,
\item $0$ above the main diagonal,
\item derivatives of $B$ below the diagonal,
\end{enumerate} we conclude that the only
part of matrix we need to take care of is
the set of first $n$ columns. Let $1 \Le q \Le n$ and
$n\cdot m+1 \Le p \Le n\cdot(m+1)$. We have:
\begin{align*}
\Delta(c_{p,q}) &= \Delta\left(\binom{i}{m}\partial^m b_{p-m\cdot n,q}\right) = \binom{i}{m}(\partial^m\Delta(b_{p-m\cdot n,q})) = \\
&=\binom{i}{m}\partial^m\left(\sum_{l=1}^{n}b_{p-m\cdot n,l}\otimes b_{l,q}\right) =\\
&= \sum_{l=1}^{n}\sum_{r=0}^m\binom{i}{m}\binom{m}{r}(\partial^rb_{p-m\cdot n,l})\otimes(\partial^{m-r} b_{l,q}) =\\
&= \sum_{l=1}^{n}\sum_{r=0}^m\binom{i-r}{m-r}\binom{i}{r}(\partial^rb_{p-m\cdot n,l})\otimes(\partial^{m-r} b_{l,q}),
\end{align*}
because $$
\binom{i}{m}\binom{m}{r} = \frac{i!m!}{m!(i-m)!r!(m-r)!}=
\frac{(i-r)!i!}{(m-r)!(i-m)!r!(i-r)!}
$$
and this is exactly what we needed.
\end{proof}

\section{Essential properties of differential representations}\label{PropertiesSection}
In the following we will use different equivalent definitions of differential
representations of a linear differential algebraic group. Summarizing
the previous sections, we see that a differential representation $V\in \Ob(\Seq)$ of a linear differential algebraic group $G$ with Hopf algebra $A$
can be defined in the following ways:
\begin{itemize}
\item by a differential morphism $G(\U) \to \GL\left(V\right)(\U)$ (Section~\ref{Representations});
\item by a natural map of group functors $G \to \Aut(V)$ (Definition~\ref{DiffModules});
\item by a differential $A$-comodule structure on $V$ (Definition~\ref{DiffComodules}, Theorem~\ref{EquivalentDefsOfRepresentations}).
\end{itemize}
Differential representations of $G$ form a category which we denote by $\Rep_G$. The objects
$\Ob(\Rep_G)$ are the underlying vector spaces.
If $V, W \in \Ob(\Rep_G)$ and $r_V : G \to \GL(V)$
and $r_W : G \to \GL(W)$ are the corresponding 
representations then $\Hom(V,W)$ consists of
those $\K$-linear maps between $V$ and $W$ that commute with the action of $G$.

\subsection{Recovering representations}\label{RecoveringSection}
Let $G$ be a linear differential algebraic group
with $A := \K\{G\}$ and $G \to \Aut(V)$ be its faithful
representation. The $\partial$-coalgebra structure on $A$ makes $A$ a $\partial$-$A$ comodule:
\begin{align}\label{RegularAction}
\rho_A := \Delta : A\to A\otimes A
\end{align}
called the {\it regular representation} of $G$.

\begin{lemma}\label{Regular} Every finite dimensional differential representation
$r_U : G \to \GL(U)$ embeds in a finite sum of copies of the regular representation of $G$.
\end{lemma}
\begin{proof}
Denote $M = U\otimes A$. Then $M$ is a differential comodule
with $$\id_U\otimes\Delta : M \to M\otimes A.$$
Since $(\id\otimes\Delta)\circ\rho = (\rho\otimes\id)\circ\rho$,
the map $\rho : U \to M$ is a map of $A$-comodules. It is injective,
because $v = (\id\otimes\varepsilon)\circ\rho(v).$ Finally,
$M \cong A^{\dim\left(U\right)}$.
\end{proof}

\begin{proposition}\label{GenerateRepresentations}
Every differential representation $U$ of $G$ is a subquotient of several copies of
a $G$-module
$$
V^{(i_1)}\otimes\ldots\otimes V^{(i_k)}\otimes V^*\otimes\ldots\otimes V^*.
$$
\end{proposition}
\begin{proof}
Fix a basis $\{u_1,\ldots,u_m\}$ of $U.$
By Lemma~\ref{Regular}
the representation $U$ is an $A := \K\{G\}$-subcomodule of
$$U\otimes A = (u_1\otimes A)\oplus\ldots\oplus (u_m\otimes A) \cong A^m.$$
Consider the canonical projections $\pi_i : A^m \to A$,
which are $G$-equivariant maps with respect to the comultiplication $\Delta : A\to A\otimes A$.
Since $U \subset A^m$, we have
$$
U \subset \bigoplus_{i=1}^m \pi_i(U)
$$
and each $\pi_i(U)$ is a $G$-module, because $\pi_i$ is
$G$-equivariant.

Consider the following surjection
$$\pi: B := \K\{X_{11},\ldots,X_{nn},1/\det\} \to A \to 0.$$
Since $\pi_i\left(U\right)$ is a finite dimensional $G$-subspace of $A$,
there exist numbers $r, s, p \in \mathbb{Z}_{\Ge 0}$ such that $\pi_i\left(U\right)$ is contained
in $\pi(L_{r,s,p}),$ where
$$L_{r,s,p} := (1/\det)^r\{f(X_{ij})\:|\:\deg(f)\Le s, \ord(f) \Le p\}.$$
There is a $B$-comodule structure on $B$ given by
\begin{align*}
\Delta(X_{ij}) &= \sum_{l=1}^nX_{il}\otimes X_{lj},\\
\Delta(\partial X_{ij}) &= \sum_{l=1}^n((\partial X_{il})\otimes X_{lj}+X_{il}\otimes (\partial X_{lj}))
\end{align*} and $L_{r,s,p}$ is a $B$-subcomodule of $B$, because of
$$
\Delta(X_{ij}X_{pq}) = \sum_{l,r = 1}^nX_{il}X_{pr}\otimes X_{lj}X_{rq}
$$
and Lemma~\ref{IthRestriction}.
We then have that $L_{r,s,p}$
is also an $A$-subcomodule of $B$. Hence, each $\pi_i(U)$ is a
subquotient of some $L_{r,s,p}$. Thus, we only need to show how
to construct these $L_{r,s,p}.$

Fix a basis $\{v_1,\ldots,v_n\}$ of $V$. We have a $B$-comodule $V$
with respect to
$$
\rho(v_j) =\sum_{i=1}^nv_i\otimes X_{ij}.
$$
For each $i,$ $1 \Le i \Le n$, the map $\varphi_i : v_j \mapsto X_{ij}$ is
$\GL_n$ (hence, $G$)-equivariant, because
$$
\varphi_i(\rho_V(v_j)) = \varphi_i\left(\sum_{l=1}^nv_l\otimes X_{lj}\right)
= \sum_{l=1}^nX_{il}\otimes X_{lj} = \Delta(X_{ij})
= \rho_B(\varphi_i(v_j))
$$ and both $\rho$ and $\Delta$ preserve the product rule with respect to $\partial $.

Consider the space of linear polynomials $L_{0,1,p}$ in the variables $\{X_{ij}\}$
and their derivatives of order up to $p$. An element $f$ of such a space
is of the form
$$
f = \sum_{i,j=1}^n\sum_{q=0}^pc_{ij}X_{ij}^{(q)},
$$
where $c_{ij} \in \K$. As it has been noticed above this space is an
$A$-subcomodule of $B$. The map $(\varphi_1,\ldots,\varphi_n)$ gives
an $A$-comodule isomorphisms between the $n$th power
$\left(V^{(p)}\right)^n$ of the $p$th derivative of the original representation
of $G$ and $L_{0,1,p}.$ Hence, one can construct
$L_{0,1,p}.$

Consider any $s \in \mathbb{Z}_{\Ge 2}$. The $G$-space $L_{0,s,p}$ is
the quotient of $\left(L_{0,1,p}\right)^{\otimes s}$ by the symmetric relations. So, we have
all $L_{0,s,p}$. Let now $s = n = \dim_{\K} V.$ Then the one-dimensional
representation $\det : G \to \K$ with $g \mapsto \det(g)$ is in $L_{0,n,p}.$
For $f \in \K^*$ we have
$$
\det(g)(f)(x) = f(x/\det(g))=\frac{1}{\det(g)}f(x).
$$
Thus,
$$
L_{r,s,p} = ({\det}^*)^{\otimes r}\otimes L_{0,s,p},
$$
which
is what we wanted to construct.
\end{proof}

\subsection{Example} We will show how Proposition~\ref{GenerateRepresentations} works
step by step.
\begin{example}
Consider the differential representation
$\rho : \Gm \to \Ga$ by
$$\Gm \ni y \mapsto
\begin{pmatrix}
1& \frac{\partial y}{y}\\
0& 1
\end{pmatrix} \in \Ga.$$
So, the underlying vector space $U$ is $\K^2$.
The representation $\rho$ corresponds to the map
of $\partial$-$\K$-algebras $$\rho^* : B := \K\left\{X_{11},X_{12},X_{21},X_{22},1/\det\right\}
\to A := \K\left\{y,1/y\right\} \to 0$$
with
\begin{align*}
X_{11} \mapsto 1, & X_{12} \mapsto y'/y, \\
X_{21} \mapsto 0, & X_{22} \mapsto 1.
\end{align*}
Take the standard basis $\{u_1,u_2\}$ of $\K^2$. We then have $\rho : U\to U\otimes A$ given by
\begin{align*}
\rho(u_1) &= u_1\otimes 1,\\
\rho(u_2) &= u_1\otimes (y'/y) + u_2\otimes 1.
\end{align*}
So, as an $A$-comodule
$$
\K^2 \subset \Span_{\K}\left\{u_1\otimes 1, u_1\otimes (y'/y)\right\}\oplus\Span_{\K}\left\{\{u_2\otimes 1\right\}.
$$
Hence, it is sufficient to construct
$\Span_{\K}\left\{u_1\otimes 1, u_1\otimes (y'/y)\right\}$ or, equivalently,
$$
W := \Span_{\K}\left\{1, \frac{y'}{y}\right\}.
$$
Consider the $B$-subcomodule $L_{0,1,0}$
of linear polynomials in $X_{11},X_{12},X_{21},X_{22}$ with
coefficients in $\K$ which is also an $A$-subcomodule of $B$.
The $A$-comodule $W$ is contained in the image of $L_{0,1,1}$ with respect to $\rho^*$.
Hence, $W$ is a subquotient of $L_{0,1,1}$.
The $A$-comodule $L_{0,1,1}$ is constructed as follows. It is enough to get
$L_{0,1,0},$ as $L_{0,1,1} = {\left(L_{0,1,0}\right)}^{(1)}.$

The group $\Gm$ has a representation on $V = \Span_{\K}\{v_1,v_2\}$
as
\begin{align*}
v_1 &\mapsto v_1\otimes y,\\
v_2 &\mapsto v_2\otimes y.
\end{align*}
Consider the two $A$-comodule maps
\begin{align*}
\varphi_1: v_1 \mapsto X_{11},\ v_2 \mapsto X_{12},\\
\varphi_2: v_1 \mapsto X_{21},\ v_2 \mapsto X_{22}.
\end{align*}
The map $(\varphi_1,\varphi_2)$ provides an $A$-comodule
isomorphism between $V^2$ and $L_{0,1,0}$. Summarizing,
we need to take the $4$th power of the original faithful representation
of $\Gm$ on $\K$, compute several subquotients, and then
sum up ($\oplus$) the result, take a subrepresentation,
and differentiate to obtain the representation $\rho$.
\end{example}

\section{Linear groups of constant elements}\label{ConstantMatricesSection}
Let $\K$ be a differential field of characteristic zero with field of constants $C$. We say that $H$ is a {\it group of constant matrices} if
it is a subgroup of some $\GL_n(C).$

\begin{proposition}\label{GroupsOfConstants} A linear differential algebraic group $G \subset \GL(V)$ is conjugate
to a group $H \subset \GL(V)(C)$ of constant matrices iff
$$
\bar{\K}\otimes V^{(p)} = \bar{\K}\otimes
\left(V\oplus\bigoplus_{i=1}^p V_i\right)
$$
for all
$p \Ge 1,$ where $V$ is a faithful representation of $G,$ $V_i \cong V,$
and $\bar{\K}$ is the differential closure of $\K$
(see, for instance, \cite[page 120]{PhyllisMichael}).
\end{proposition}
\begin{proof}
Let $\{v_1,\ldots, v_n\}$ be a $\K$-basis of $V.$ Assume that there exists
a matrix $D \in \GL_n\left(\bar{\K}\right)$ such that
$$
D^{-1} GD = H.
$$
For $g \in G$ let $A_g$ be the corresponding matrix with respect to
the basis $\{v_1,\ldots, v_n\}.$ The matrix of $g$ with respect to
the basis
$$
(w_1,\ldots,w_n) = (v_1,\ldots,v_n)\cdot D
$$
is given by
$$
B_g = D^{-1} A_gD.
$$
Hence, $B_g \in \GL_n(C).$
We have
$$
g\cdot (\partial^p\otimes w_i) = \partial^p\otimes (g\cdot w_i) = \partial\otimes(B_g w_i) = B_g(\partial\otimes w_i).
$$
Thus,
\begin{align*}
\bar{\K}\otimes V^{(p)} &= \Span_{\bar{\K}}\{w_1,\ldots,w_n\}\oplus
\ldots\oplus \Span_{\bar{\K}}\{\partial^p\otimes
w_1,\ldots,\partial^p\otimes w_n\} = \\
&=\bar{\K}\otimes
\left(V\oplus\bigoplus_{i=1}^p V_i\right),
\end{align*}
where $V_i = \Span_{\K}\{\partial^i\otimes
w_1,\ldots,\partial^i\otimes w_n\}.$

Conversely, let
$$
\bar{\K}\otimes V^{(1)} = (\bar\K\otimes V)\oplus V_1,
$$
where $V_1 \cong \bar{\K}\otimes V.$
Choose bases $\{v_1,\ldots,v_n\}$ and $\{w_1,\ldots,w_n\}$ in $V$ and $V_1,$ respectively.
There exists a matrix $A \in \Mn_{(2\cdot n)\times n}(\K)$ such that
$$
\left(\partial\otimes v_1,\ldots,\partial\otimes v_n\right) =
\left(v_1,\ldots,v_n,w_1,\ldots,w_n\right)\cdot A.
$$
For a matrix $B \in \GL_n\left(\bar{\K}\right)$ we have
\begin{align*}
\partial\otimes((v_1,\ldots,v_n)\cdot B) &= \left(\partial\otimes (v_1,\ldots,v_n)\right)\cdot B + \left(v_1,\ldots,v_n\right)\cdot\partial(B) =\\
&=\left(v_1,\ldots,v_n,w_1,\ldots,w_n\right)\cdot AB + \left(v_1,\ldots,v_n\right)\cdot\partial(B).
\end{align*}
We will first show that there are $n\times n$ invertible matrices $B$ and $D$ such that
\begin{equation}\label{eq1}
\partial\otimes((v_1,\ldots,v_n)\cdot B) = (w_1,\ldots,w_n)\cdot D.
\end{equation}
It follows then that
$$
AB + \begin{pmatrix}\partial B\\0\end{pmatrix} = \begin{pmatrix}
0\\
D
\end{pmatrix}.
$$
Let $A = \begin{pmatrix} A_1\\A_2\end{pmatrix}.$ We then obtain the following
system:
\begin{equation}\label{eq3}
\begin{cases}
\partial B &= -A_1B;\\
A_2 B  &= D.
\end{cases}
\end{equation}
Since the field $\bar\K$ is differentially closed, there exists a matrix
$B \in \GL_n\left(\bar{\K}\right)$ satisfying the first equation of system~\eqref{eq3}.
Now, equation~\eqref{eq1} forces $D$ to be invertible
as well, as the dimension of the span of $\partial\otimes((v_1,\ldots,v_n)\cdot B)$ is equal to $n.$

Let
$$
(u_1,\ldots,u_n) = (v_1,\ldots,v_n)\cdot B.
$$
Consider $g \in G$ and its matrix $A_g =(a_{ij})$ with respect to the basis $\{u_1,\ldots,u_n\}.$ For each $j,$ $1\Le j \Le n,$ we have:
$$
g\cdot\partial\otimes u_j =\partial\otimes\left(\sum_{i=1}^n a_{ij}u_i\right) =
\sum_{i=1}^n\left(\partial(a_{ij})\otimes u_i + a_{ij}\cdot\partial\otimes u_i\right).
$$
By the construction,
the space
$$
\Span_{\K}\left\{\partial\otimes u_1,\ldots,\partial\otimes u_n\right\}
$$
is $G$ invariant. This implies that $\partial(a_{ij}) = 0$
for all $i$ and $j,$ $1 \Le i,j \Le n.$ Thus, $A_g \in \Mn_n(C).$
\end{proof}

\begin{corollary} A linear differential algebraic group $G \subset \GL(V)$ is
differentially isomorphic over $\bar{\K}$ to an algebraic subgroup of some
$\GL_m(C)$ if and only if {\bf there exists} a faithful
representation $W$ of $G$ such that
$$
\bar{\K}\otimes W^{(1)} = \bar{\K}\otimes(W\oplus W_1),
$$
where $W_1 \cong W$ as differential representations of the group $G$.
\end{corollary}
\begin{proof}
Let $\bar{\K}\otimes W^{(1)} = \bar{\K}\otimes (W\oplus W_1)$. The representation morphism  $r_W : G
\to \GL(W)$ is a differential algebraic group homomorphism and is a
differential isomorphism between $G$ and the image $r_W(G).$ By
Proposition~\ref{GroupsOfConstants} the group $r_W(G)$ is conjugate
to a group of matrices with constant coefficients. Composition of
this isomorphism with $r_W$ gives the desired differential algebraic
group isomorphism.

Let now $G$ be differentially isomorphic to a subgroup of $\GL_m(C)$.
This gives a faithful representation $W$ of $G$. Moreover, since the
matrices have constant entries, we have $\bar{\K}\otimes W^{(1)} = \bar{\K}\otimes (W\oplus W)$.
\end{proof}

\section{Tannaka's theorem for linear differential algebraic groups}\label{TannakasTheoremSection}
In this section we will show how one can recover a linear
differential algebraic group knowing all its representations. For this
we first prove Tannaka's theorem (Theorem~\ref{TannakaTheorem}). Then
using this fact we reconstruct the differential Hopf algebra of functions on the
group in Sections~\ref{RecoverSection} and \ref{RecoverDeltaS}.
Parts of the proof closely follow \cite[Section 2.5]{Springer}.
The novelty lies in the fact that we can recover the differential
structure on the Hopf algebra.

\subsection{Preliminaries}

Let $G$ be a linear differential algebraic group with the Hopf algebra $A := \K\{G\}$. Note that $A$ is also a locally finite $G$-module with the action~\eqref{RegularAction} by \cite[Theorem, page 230]{CassidyRep}.
Let $\omega : \Rep_G \to \Seq$
be the forgetful functor.

\begin{definition} For a $\partial$-$\K$-algebra $B$
we define the group $\AutOP(\omega)(B)$ to be the set of sequences
$$\lambda(B) = (\lambda_X\:|\: X\in \Ob(\Rep_G)) \in
\AutOP(\omega)(B)$$ such that $\lambda_X$ is a $B$-linear automorphism
of $\omega(X)\otimes B$ for each $G$-space $X$, $\omega(X) \in \Ob(\Seq)$, that is,
$\lambda_X \in \Aut_B(\omega(X)\otimes B)$, such that
\begin{itemize}
\item for all $X_1,$ $X_2$ we have
\begin{align}\label{TensorSpreading}
\lambda_{X_1\otimes X_2} = \lambda_{X_1}\otimes\lambda_{X_2},
\end{align}
\item $\lambda_{\underline{1}}$ is the identity map on
$\underline{1}\otimes B = B$,
\item for every $\alpha \in \Hom_G(X,Y)$ we have
\begin{align}\label{Equivariance}
\lambda_Y\circ(\alpha\otimes \id_B) = (\alpha\otimes \id_B)\circ\lambda_X :
X\otimes B \to Y\otimes B,
\end{align}
\item for every $X$ we have
\begin{align}\label{CommuteWithD}
\partial\circ\lambda_X = \lambda_{X^{(1)}}\circ \partial,
\end{align}
\item the group operation $\lambda_1(B)\cdot\lambda_2(B)$ is defined
by composition in each set $\Aut_B(\omega(X)\otimes B).$
\end{itemize}
\end{definition}
$\AutOP(\omega)$ is a functor from the category
of $\partial$-$\K$-algebras to groups:
$$
B \mapsto \AutOP(\omega)(B).
$$ 
Any $g \in G(B)$ determines an element $\lambda_g \in \AutOP(\omega)(B)$,
because if $X \in \Ob(\Rep_G)$ and $\Phi_B : G(B) \to \GL(\omega(X)\otimes B)$ then $\Phi(g)$ is a
$B$-linear automorphism of $\omega(X)\otimes B$ and the property~\eqref{Equivariance}
is satisfied by the definition of $G$-equivariance of $\alpha$. So,
we have a morphism $\Phi$ of functors
$$G \to \AutOP(\omega), \quad g\in G(B) \mapsto \Phi(g) \in \AutOP(\omega)(B)$$
for any $\partial$-$\K$-algebra $B$ as for any $\varphi : B_1 \to B_2$
we have the following commutative diagram:
\begin{equation*}
\begin{CD}
G(B_1) @>{\Phi_{B_1}}>>\AutOP(\omega)(B_1) \\
@VV{G(\varphi)}V @VV{\AutOP(\omega)(\varphi)}V \\
G(B_2) @>{\Phi_{B_2}}>> \AutOP(\omega)(B_2)
\end{CD}\qquad\qquad
\end{equation*}
where $G(\varphi)$ and $\AutOP(\omega)(\varphi)$ denote the
morphisms
$$
[\psi : A \to B_1] \mapsto [G(\varphi)(\psi)= \varphi\circ\psi: A \to B_2 ];
$$
\begin{align*}
\left[\lambda(B_1)\right. &=\left. \left(\lambda_X(B_1) : \omega(X)\otimes B_1 \to \omega(X)\otimes B_1\right)\right] \mapsto \\
\left[\AutOP(\omega)(\varphi)(\lambda(B_1))\right. &= \lambda(B_2)=\\
&= \left.\left(\left(\id_{\omega(X)}\otimes\varphi\right)\circ\lambda_X(B_1) : \omega(X)\otimes B_2 \to \omega(X)\otimes B_2\right)\right],
\end{align*}
respectively. The latter means that we take the restriction of
$\lambda_X(B_1)$ to $\omega(X)$ and map it to $\omega(X)\otimes B_1$ and then apply
$\id_{\omega(X)}\otimes\varphi$ mapping
it to $\omega(X)\otimes B_2$. At the end we prolong such a
map to $\omega(X)\otimes B_2$ by $B_2$-linearity.

\subsection{The theorem}
\begin{theorem}\label{TannakaTheorem} For a linear differential algebraic group $G$ let $\omega : \Rep_G \to \Seq$  be the forgetful
functor. Then
$$ G \cong \AutOP(\omega)
$$ as functors.
\end{theorem}
\begin{proof} (Following \cite[Theorem 2.5.3]{Springer} with
{\it differential} modification)
Since $g \in G$ determines an element of $\AutOP(\omega)$, we
only need to show the converse. Let $B$ be a
$\partial$-$\K$-algebra and
$(\lambda_X) \in \AutOP(\omega)(B)$.
Let $V \in \Ob(\Rep_G)$  not necessarily finite dimensional
but locally finite. This means that every vector $v \in V$
is contained in a differential $G$-module $W$ with $\dim_{\K}W < \infty.$

We will show that for given $B,$ $\lambda,$  and a locally
finite differential $G$-module $V$, there exists a $B$-linear automorphism of $V\otimes B$ (denoted by $\lambda_V$) such that
the properties \eqref{TensorSpreading}, \eqref{Equivariance}, and \eqref{CommuteWithD} are satisfied
and for any $W \subset V$ we have ${\lambda_V|}_W = \lambda_W.$
For $v \in V$ we take $W \in \Ob(\Rep_G)$ such that
$v \in W$ and $\dim_{\K} W < \infty$.

First, define
$\lambda_V(v) := \lambda_W(v)$. We need to show the correctness.
Let $W'$ be another representation such that $v \in W'$. Consider the $G$-module
$W\cap W' \ni Gv$. From \eqref{Equivariance}
it follows that 
$$
\lambda_{W'}(v) = \lambda_{W\cap W'}(v) = \lambda_W(v).
$$
Since each $\lambda_W$ is invertible and linear, the map $\lambda_V$
is also invertible and linear.
We are going to show now that \eqref{TensorSpreading},  \eqref{Equivariance}, and \eqref{CommuteWithD} hold for locally finite modules.
For locally finite $V$ and $V'$
we choose $W\ni v$ and $W' \ni v',$ objects of $\Seq$. Then
$$
\lambda_{V\otimes V'}(v\otimes v') = \lambda_{W\otimes W'}(v\otimes v')
= \lambda_W(v)\otimes\lambda_{W'}(v') = \lambda_V(v)\otimes\lambda_{V'}(v').
$$
Hence, the property \eqref{TensorSpreading} is
satisfied. Also, let $\alpha \in \Hom_G(V,V').$ Consider $v \in V$
and $W$ such that $v \in W$ and $\dim W < \infty$ and of the same
kind $W' \supset \alpha(W).$ We obtain that
$$
\lambda_{V'}\circ\alpha(v) = \lambda_{W'}\left({\alpha\big|}_W(v)\right) ={\alpha\big|}_W\circ\lambda_W(v)
=\alpha\circ\lambda_V(v).
$$
We then have \eqref{Equivariance}. Moreover, since $W^{(1)} \subset V^{(1)},$ we have
$$
\partial\circ\lambda_V (v) = \partial\circ\lambda_W(v) = \lambda_{W^{(1)}}(\partial v) = \lambda_{V^{(1)}}(\partial v),
$$
which implies \eqref{CommuteWithD}.
Thus, we may say $\lambda_V \in \Aut_B(V\otimes B)$.

Recall that the $G$-module $A$ is locally finite via
$$
\rho : A\times G \to A, \quad (\rho(g)f)(x) = f(x\cdot g)
$$
by
\cite[Theorem, page 230]{CassidyRep}, where $x, g \in G(B)$ and $f \in A = \K\{G\}.$
The same is true for $A\otimes A$. Consider the multiplication map
$$m : A\otimes A \to A, \quad f\otimes h \mapsto f\cdot h$$ which is
$G$-equivariant, because each $g \in G(B)$
is a ($B$-linear) algebra automorphism of $A\otimes B$. According to
\eqref{TensorSpreading} and \eqref{Equivariance} we have
$$
 m\circ (\lambda_A\otimes \lambda_A) = m\circ\lambda_{A\otimes A}
= \lambda_A\circ m.
$$
Moreover, from \eqref{CommuteWithD} we conclude that
$\lambda_A$ is a
$\partial$-$\K$-algebra automorphism
$A \to A$.
We will show that this must correspond to an element in the group. More precisely,
there exists a differential algebraic automorphism
$\varphi : G(B)\to G(B)$ such that for all $f \in A$ and $g \in G(B)$ we have
$$\lambda_A(f)(g) = f(\varphi(g)).$$
We will show that this morphism $\varphi$ is right
multiplication by an element of $G(B)$. This will show that
an element of the group corresponds to $\lambda_A$. After that
we demonstrate that the algebra $A$ can be replaced by any $G$-module.

For every $f \in A$ and $g,$ $h \in G(B)$ we have
$\Delta(f)(g,h) = f(gh).$ Take any $f \in A$ and $g, x, y \in G(B)$. We
have
$$\Delta\circ\rho(g)(f)(x,y)= f(xyg) =
(\id_A\otimes\rho(g))\circ\Delta(f)(x,y).$$
Hence,
\begin{align}\label{eqDelta}
\Delta\circ\rho(g) = (\id_A\otimes\rho(g))\circ\Delta.
\end{align}
Consider the locally finite $G$-module $U := A\otimes A$ via
$r_U = \id_A\otimes\rho.$ Then, by~\eqref{eqDelta} the map $\Delta$
is $G$-equivariant for $\rho$ and $r_U$. From~\eqref{Equivariance}
we have $\Delta\circ\lambda_A = \lambda_U\circ\Delta.$
Because of $\eqref{TensorSpreading}$ we obtain that
$$
\lambda_U = \lambda_{A, \id_A}\otimes\lambda_{A,\rho} = \id\otimes\lambda_{A,\rho},
$$
because $\lambda_I = \id$ and $\lambda$ is $\K$-linear. Thus,
$$
\Delta\circ\lambda_A = (\id\otimes\lambda_A)\circ\Delta.
$$
For any $f \in A$ and $g,$ $h\in G$ we have $\Delta\circ\lambda_A(f)(g,h)=
f(\varphi(gh)).$ On the other hand,
$$
(\id\otimes\lambda_A)\circ\Delta(f)(g,h) = f(g\varphi(h)).
$$
Thus,
$$
\varphi(gh) = g\varphi(h).
$$
Let
$$
x = \varphi(e),
$$
which is a {\it differential} algebra homomorphism $A \to B.$
From
this we conclude that for any $g \in G(B)$ one has $\varphi(g) = gx$.
Hence, $\lambda_A = \rho(x)$. It remains to show that other
automorphisms $\lambda_V$ look the same
(completely determined by this element $x$).

Consider any $V \in \Ob(\Rep_G)$ with the action $r_V$. For any $u \in V^*$ there is a $G$-homomorphism
$$\varphi_u : V \to A,\quad v \mapsto \varphi_u(v), \quad \varphi_u(v)(g) := u(\rho(g)\cdot v),$$
where
$v \in V.$
By~\eqref{Equivariance} and the above, we have $\rho(x)\circ\varphi_u = \lambda_A\circ\varphi_u = \varphi_u\circ\lambda_V.$
Take any $g \in G(B)$ and $v \in V$. We then have
\begin{align*}
\rho(x)(u(r_V(g)(v))) = u(r_V(gx)(v)),\\
\varphi_u\circ\lambda_V(v)(g)= u(r_V(g)\circ\lambda_V(v)).
\end{align*}
Thus,
$$
r_V(x) = r_V(ex) = r_V(e)\circ\lambda_V = \lambda_V,
$$
because the elements $v,$ $g,$ and $u$ were arbitrary.
\end{proof}

\subsection{Recovering the differential Hopf algebra of $G$}\label{RecoverSection}
Similar to \cite[Sections 2.5.4--2.5.8]{Springer} we can recover the Hopf algebra $A = \K\{G\}$ in the
following way. In addition, we show how to obtain the {\it differential
structure} on $A$ (see Lemma~\ref{DiffStructureLemma}).

\subsubsection{First step} We will construct the map $\psi_V$
from ``some representations'' of $G$ to the algebra $A$ of regular
differential functions on $G$ and study main properties of $\psi_V.$

Recall that for $V \in \Seq$ we denote $V^{(i)} = \K[\partial]_{\Le i}\otimes V$ (non-commutative tensor product) and sometimes we
write $V^{(0)}$ instead of $V$ for convenience.
For $V \in \Ob(\Rep_G)$ and
\begin{align*}
v \in V,\ u \in V^*
\end{align*}
we have the linear map
\begin{align}\label{SpringersMap}
\psi_V : \ V\otimes V^* \to A, \quad \psi_V(v\otimes u)(g) = u(r_V(g)\cdot v).
\end{align}
Also we introduce the following map:
$$
F: V^* \to \left(V^{(1)}\right)^*,\quad F(u)(v) = u(v),\ F(u)(\partial\otimes v) = \partial(u(v)),\ v \in V.
$$

\begin{lemma}\label{PhiLemma}We have the following properties:
\begin{enumerate}
\item If $\phi \in \Hom_G(V,W)$ then
$$
\psi_V\circ(\id\otimes \phi^*) = \psi_W\circ(\phi\otimes \id)
$$
as maps of $V\otimes W^* \to A.$
\item We have
$$
\psi_{V\otimes W} = m\circ(\psi_V\otimes\psi_W)\circ c,
$$
where $c : (V\otimes W)\otimes (V\otimes W)^* \cong (V\otimes V^*)\otimes(W\otimes W^*).$
\item Moreover,
$$
\partial(\psi_V(v\otimes u)) = \psi_{V^{(1)}}\left((\partial v)\otimes F(u) \right).
$$
\end{enumerate}
\end{lemma}
\begin{proof}
For $v \in V,$ $u \in V^*,$ and $g \in G$ we have
\begin{align*}
\psi_V\circ(\id\otimes\phi^*)(v\otimes u)(g) &= \psi_V(v\otimes\phi^*(u))(g) =
\phi^*(u)(r_V(g)\cdot v) = \\
&=u(\phi(r_V(g)\cdot v)) = u(r_W(g)\cdot\phi(v)) =\\
&=\psi_W\circ(\phi\otimes\id)(v\otimes u)(g),
\end{align*}
Furthermore, consider $w \in W$ and $t \in W^*$. We then also have
\begin{align*}
\psi_{V\otimes W}(v\otimes w\otimes u\otimes t)(g) &= (u\otimes t)(r_{V\otimes W}(g)\cdot(v\otimes w)) =\\
&= (u\otimes t)((r_V(g)\cdot v)\otimes(r_W(g)\cdot w)) =\\
&= u(r_V(g)\cdot v)\cdot t(r_W(g)\cdot w) =\\
&= m\circ(\psi_V\otimes\psi_W)((v\otimes u)\otimes(w\otimes t))(g))=\\
&= m\circ(\psi_V\otimes\psi_W)\circ c(v\otimes w\otimes u\otimes t)(g).
\end{align*}
Let $\{e_1,\ldots,e_n\}$ be a basis of $V$ with the dual basis $\{f_1,\ldots,f_n\}$ and take $f \in \K.$
We have:
\begin{align*}
\partial(\psi_V(f\cdot e_i\otimes f_j))(g) &= \partial(f\cdot f_j(r_V(g)\cdot e_i))= \partial\left(f\cdot g_{ij}^V\right) =\\
&=\partial(f)\cdot g_{ij}^V +  f\cdot\partial\left(g_{ij}^V\right)=\\
&= \partial(f)\cdot f_j(r_V(g)\cdot e_i) + f\cdot\partial(f_j(r_V(g)e_i)) =\\
&=\partial(f)\cdot\psi_{V^{(1)}}(e_i\otimes f_j)(g)+f\cdot\psi_{V^{(1)}}((\partial e_i)\otimes F(f_j))(g)=\\
&= \psi_{V^{(1)}}\left(\partial(f\cdot e_i)\otimes F(f_j)\right)(g).
\end{align*}
\end{proof}

\subsubsection{Second step} Here, we will construct a differential algebra $\A$ (that will be our candidate for $A$) using representations of $G$ as objects of $\Seq$ together with
morphisms between them and not using any other information from $G$.

Let
$$
\F = \bigoplus_{V\in \Ob(\Rep_G)} V\otimes V^*.
$$
So, the canonical injections
$$
i_V : V\otimes V^*\to \F
$$
are defined.
Consider the subspace $\Rr$ of $\F$ spanned by
\begin{align*}
\left\{\left(i_V(\id\otimes \phi^*) - i_W(\phi\otimes
\id)\right)(z)\:\big|\: V,\:W \in \Ob(\Rep_G),\phi\in
\Hom\left(V,W\right), z \in V\otimes W^*\right\}.
\end{align*}
We now put $\A = \F/\Rr.$ For $v  \in V$ and $u \in V^*$ we denote by
$$
a_V\left(v\otimes u\right)
$$
the image in $\A$ of
$$
i_V\left(v\otimes u\right).
$$
So, for any $\phi \in \Hom_G(V,W)$ we have
\begin{align}\label{undermorphisms}
a_V(v\otimes \phi^*(u)) =
a_W(\phi(v)\otimes u).
\end{align}
\begin{lemma} For all $v \in V$ and $u \in V^*$ we have
\begin{equation}\label{NewDiffEq}
a_V(v\otimes u)= a_{V^{(1)}}(\partial v\otimes \varphi^*(u)),
\end{equation}
where the morphism
$$
\varphi : V^{(1)}\to V,\quad 1\otimes v \mapsto 0,\ \partial\otimes
v \mapsto v.
$$
\end{lemma}
\begin{proof}
Follows from formula~\eqref{undermorphisms}.
\end{proof}

Take $v \in V,$ $w \in W,$ $u \in V^*,$
$t \in W^*.$ Let also $\{v_i\}$ be a basis of $V$ and $\{u_i\}$ be its dual.
Introduce the following operations on $\A:$
\begin{align}
m(a_V(v\otimes u), a_W(w\otimes t)) &= a_{V\otimes W}((v\otimes w)\otimes(u\otimes t)),\label{MultFormula}\\
\partial(a_V(v\otimes u)) &= a_{V^{(1)}}(\partial v\otimes F(u)),\label{DiffFormula}\\
\tilde{\Delta}(a_V(v\otimes u)) &= \sum_j a_V(v_j\otimes u)\otimes a_V(v\otimes u_j),\label{DeltaFormula}\\
\tilde{S}(a_V(v\otimes u)) &= a_{V^*}(u\otimes v).\label{SFormula}
\end{align}

\begin{proposition} The $\K$-vector space $\A$ contains a non-zero vector.
\end{proposition}
\begin{proof} Let $\Cat$ be the category of differential
representations of the trivial differential group $G_e = \{e\}.$ Note
that the algebra of $G_e$ is just the differential field $\K.$ Let
$V,\:W \in \Ob(\Cat)$. Choose ordered bases $\{v_1,\ldots,v_n\}$ and
$\{w_1,\ldots,w_n\}$ of $V$ and $W$, respectively. Take
$a_V(v_1\otimes v_2^*)$ and $a_W(0\otimes w_2^*)$. Consider the
linear map:
$$\phi: V \to W,\quad v_2\mapsto w_2,\ v_i \mapsto 0,\ i\ne 2.$$
We have $\phi^*(w_2^*) = v_2^*.$ Then,
$$
a_V(v_1\otimes v_2^*) = a_V(v_1\otimes\phi^*(w_2^*)) = a_W(0\otimes w_2^*) = 0.
$$
Let now $0\ne v \in V$ and $0\ne w \in W.$ Without a loss of
generality we may assume that $v=v_1$ and $w = w_1$. Consider the
linear map
$$
\phi : V\to W,\quad v_1\mapsto w_1,\ v_i\mapsto 0,\ i\ne 1.
$$
We have $\phi^*(w^*) = v^*.$ Then,
$$
a_V(v\otimes v^*) = a_V\left(v\otimes\phi^*(w^*)\right)= a_W(w\otimes w^*).
$$
Thus, if for the trivial representation $\1 = \Span_{\K}\{e\}$ then
$$
\A_{G_e} = \Span_{\K}\left\{a_{\1}(e\otimes e^*)\right\}.
$$
But $a_{\1}(e\otimes e^*) \ne 0$ in $\A_{G_e}.$
The vector space $\A_{G_e}$ is just the quotient of $\F$ by the
subspace generated by the  relations coming from all possible
morphisms $\phi: V\to W.$ And the latter subspace is not the whole
$\F$ as it has been shown above. Since there is a surjective linear
map $\A \to \A_{G_e},$ this implies that $\A$ contains a non-zero
vector.
\end{proof}

\begin{lemma}
The definition of $m$ is correct and provides a structure of
a commutative associative algebra on the $\K$-vector space $\A$ and
the unit $1$ is given by the trivial representation $\1.$
\end{lemma}
\begin{proof}
Consider morphisms $\phi_1 : V \to X$ and $\phi_2 : W \to Y$
as $G$-vector spaces
and vectors
$$
v \in V,\ u \in V^*,\ w \in W,\ t \in W^*,\
x \in X^*,\ y \in Y^*
$$
such that $\phi_1^*(x) = u$ and $\phi_2^*(y) = t.$ We then have:
\begin{align*}
a_V(v\otimes u)\cdot a_W(w\otimes t) =
&a_V(v\otimes\phi_1^*(x))\cdot a_W(w\otimes \phi_2^*(y))=\\
&= a_X(\phi_1(v)\otimes x)\cdot a_Y(\phi_2(w)\otimes y)=\\
&= a_{X\otimes Y}((\phi_1(v)\otimes\phi_2(w))\otimes(x\otimes y))=\\
&= a_{X\otimes Y}((\phi_1\otimes\phi_2)(v\otimes w)\otimes(x\otimes y))=\\
&=a_{V\otimes W}((v\otimes w)\otimes(u\otimes t)).
\end{align*}
We now prove that the multiplication is associative and commutative.
Consider the morphism
$$\phi \in \Hom(V\otimes W,W\otimes V),\quad v\otimes w \mapsto w\otimes v.$$
We then have
\begin{align*}
a_V(v\otimes u)\cdot a_W(w\otimes t)
&= a_{V\otimes W}((v\otimes w)\otimes(u\otimes t))=
a_{W\otimes V}((w\otimes v)\otimes (t\otimes u)) = \\
&= a_W(w\otimes t)\cdot a_V(v\otimes u).
\end{align*}
So, the multiplication is commutative.

For $X \in \Ob(\Rep_G)$ and $x \in X,$
$y\in X^*$ we also have
\begin{align*}
&(a_V(v\otimes u)\cdot
a_W(w\otimes t))\cdot
a_X(x\otimes y)
=\\
&=a_{V\otimes W}((v\otimes w)\otimes(u\otimes t))\cdot
a_X(x\otimes y)=\\
&= a_{(V\otimes W)\otimes X}
((v\otimes w)\otimes x)\otimes((t\otimes u)\otimes y)) = \\
&= a_{V\otimes (W\otimes X)}
((v\otimes (w\otimes x))\otimes(t\otimes (u\otimes y))) =\\
&= a_V(v\otimes t)\cdot a_{W\otimes X}
((w\otimes x)\otimes(u\otimes y))=\\
&=a_V(v\otimes u)\cdot(a_W(w\otimes t)\cdot
a_X(x\otimes y)).
\end{align*}
We have shown that the multiplication is associative.

Let $\1$ be the trivial representation of the group $G$ and
$0 \ne e \in \1,$ $f\in \1^*,$ $f(e) = 1.$ We have the morphism
$$
\varphi \in \Hom(V, V\otimes \1),\quad v \mapsto v\otimes e.
$$
Then
\begin{align*}
a_V(v\otimes u)\cdot a_{\1}(e\otimes f) &=
a_{V\otimes\1}((v\otimes e)\otimes(u\otimes f)) =
a_{V\otimes \1}(\varphi(v)\otimes(u\otimes f)) =\\
&= a_V(v\otimes \varphi^*(u\otimes f)) = a_V(v\otimes u).
\end{align*}
Thus, $\A$ is a commutative associative algebra with unity.
\end{proof}

Our main goal is to recover the {\it differential} Hopf algebra $A.$
We give a differential structure on $\A$ and then show that this
structure corresponds to the one of $A.$ Let us describe the {\it intuition}
behind the construction we are going to present.

We are recovering a subgroup of $\GL_n.$ Let us denote the corresponding matrix coordinate
functions by $y_{ij}.$ We must be able to:
\begin{enumerate}
\item differentiate these functions $y_{ij}$ obtaining $y_{ij}',\ldots,y_{ij}^{(p)},\ldots;$
\item multiply the results of this differentiation and stay in the algebra.
\end{enumerate}
For bases $v_1,\ldots,v_n$ and $u_1,\ldots,u_n$ of $V$ and $V^*,$ respectively, the coordinate functions mapping $G \to \K$ are given by
$$
y_{ij}(g) = \psi_V(v_i\otimes u_j)(g) = u_j(r_V(g)\cdot v_i),
$$
where $g \in G.$
The candidates for these functions in $\A$ are,
certainly, $a_V(v_i\otimes u_j).$ Our correspondence between
$\A$ and $A$ must preserve differentiation. So, $y_{ij}'$ corresponds
to $\partial(a_{V}(v_i\otimes u_j))$ which we still need to
define.

Moreover, such a definition must leave us in the same
category and satisfy the product rule for differentiation. We also
notice that since $y_{ij}',\ldots,y_{ij}^{(p)},\ldots$ and
$y_{ij}^{(q)}\cdot y_{kl}^{(r)}$ are monomials, we should preserve this property for $\partial^q\left(a_V(v_i\otimes u_j)\right)\cdot\partial^r\left(a_V(v_k\otimes u_l)\right).$ These ideas are implemented in Lemma~\ref{DiffStructureLemma} and Theorem~\ref{Recover}.

\begin{lemma}\label{DiffStructureLemma}
The natural differential structure on $\A$ introduced in formula~\eqref{DiffFormula}
makes it a $\partial$-$\K$-algebra.
\end{lemma}
\begin{proof}
First of all, recall that
$$\partial\left(a_V(v\otimes u)\right) = a_{V^{(1)}}((\partial v)\otimes F(u)).$$
Recall also that we let $F(u)(\partial v) = \partial(u(v))$ for any
$v \in V$ and $u \in V^*.$
We need to show its correctness with respect to the morphisms. Let
$t \in W^*$. We have:
\begin{align*}
\partial(a_W(\phi(v)\otimes t)) &= a_{W^{(1)}}(\partial(\phi(v))\otimes F(t)) =
a_{W^{(1)}}(\phi(\partial v)\otimes F(t)) = \\
&=a_{V^{(1)}}(\partial v \otimes \phi^*(F(t)))=
\partial(a_V(v\otimes\phi^*(t)))
\end{align*}
for a morphism $\phi : V\to W$ that we naturally prolong to a morphism
$\phi : V^{(1)} \to W^{(1)}$ mapping $\partial v \mapsto \partial(\phi(v)).$
Hence, the differentiation is correct.

We need to show the product rule. We have:
\begin{align*}
&\partial \left(a_V(v\otimes u)\cdot
a_W(w\otimes t)\right) = \partial \left(a_{V\otimes W}((v\otimes w)\otimes(u\otimes t))\right) =\\
&= a_{(V\otimes W)^{(1)}}
(\partial (v\otimes w)\otimes F(u\otimes t))=\\
&=a_{V^{(1)}\otimes W^{(1)}}
((\partial v\otimes w)\otimes (F(u)\otimes F(t))+ (v\otimes\partial w)\otimes (F(u)\otimes F(t)))=\\
&=a_{V^{(1)}}(\partial v\otimes F(u))\cdot a_W(w\otimes
t)+a_V(v\otimes u)\cdot
a_{W^{(1)}}(\partial w\otimes F(t)) =\\
&=\left(\partial a_V(v\otimes u)\right)\cdot
a_W(w\otimes t) + a_V(v\otimes u)\cdot
\left(\partial a_W(w\otimes t)\right),
\end{align*}
where we have used the morphism:
\begin{align*}
&\psi : (V\otimes W)^{(1)} \to V^{(1)}\otimes W^{(1)},\\
&\psi : 1\otimes x\otimes z \mapsto (1\otimes x)\otimes(1\otimes z),\\
&\psi: \partial\otimes x\otimes z \mapsto (\partial x)\otimes z +
x\otimes(\partial z).
\end{align*}
Its dual
$$\psi^* : \left(V^{(1)}\otimes W^{(1)}\right)^* \to \left((V\otimes W)^{(1)}\right)^*
$$
maps
$$
F(u)\otimes F(t)\mapsto F(u\otimes t)
$$
Indeed,
\begin{align*}
&\psi^*(F(u)\otimes F(t))(1\otimes x\otimes z) = (F(u)\otimes
F(t))((1\otimes x)\otimes(1\otimes z)) =\\
&= F(u)(1\otimes x)\cdot F(t)(1\otimes z) = u(x)\cdot t(z) =
F(u\otimes t)(1\otimes x\otimes z)
\end{align*}
and
\begin{align*}
&\psi^*(F(u)\otimes F(t))(\partial\otimes x\otimes z) = (F(u)\otimes
F(t))(\partial x\otimes z + x\otimes\partial z)) =\\
&= F(u)(\partial x)\cdot F(t)(1\otimes z) + F(u)(1\otimes x)\cdot
F(t)(\partial z)=\\
&= \partial(u(x))\cdot t(z) + u(x)\cdot\partial(t(z)) =
\partial(u(x)\cdot t(z)) = F(u\otimes t)(\partial\otimes x\otimes z).
\end{align*}
\end{proof}

\subsubsection{Third step} Now, we can show that the differential algebra
$\A$ we have constructed is what we were looking for.

\begin{lemma}\label{RigidLemma} Let $\Cat$ be a rigid
abelian tensor category with a tensor $\K$-linear functor $\omega: \Cat \to \Seq$
then
$$\End^\otimes(\omega) = \Aut^\otimes(\omega).$$
\end{lemma}
\begin{proof}
We expand the proof that appears in
\cite[Proposition 1.13]{Deligne}. Let $\lambda : \omega \to \omega$.
For each $X \in \Cat$ there exists a morphism $t_X : \omega(X)^* \to \omega(X)^*$ such that
the following diagram is commutative:
\begin{equation*}
\begin{CD}
\omega(X^*)@>{\lambda_{X^*}}>> \omega(X^*)\\
@V{\varphi}V{\cong}V @V{\varphi}V{\cong}V\\
\omega(X)^* @>{t_X}>>\omega(X)^*
\end{CD}
\end{equation*}
The category $\Seq$ is rigid. So, for all $U,$ $V \in \Ob(\Seq)$ we have
$\IntHom(U,V) \cong \IntHom(V^*,U^*)$. We then let
$\mu_X := (t_X)^* : \omega(X) \to \omega(X)$. For any $f : X\to Y$
the following diagram commutes:
\begin{equation*}
\begin{CD}
Y@>{}>> Y^*@>{\omega}>>\omega(Y^*)@>{\lambda_{Y^*}}>> \omega(Y^*)\\
@AA{f}A @VV{f^*}V@VV{\omega(f^*)}V @VV{\omega(f^*)}V\\
X @>{}>>X^*@>{\omega}>>\omega(X^*) @>{\lambda_{X^*}}>>\omega(X^*)
\end{CD}
\end{equation*}
Gathering all commutative diagrams together we obtain that $\mu = (\mu_X) \in \End^\otimes(\omega).$
We now
show that $\mu = \lambda^{-1}$.
We have
\begin{equation*}
\begin{CD}
X^*\otimes X@>{\omega}>> \omega(X^*\otimes X)@>{\lambda_{X^*\otimes X}}>>\omega(X^*\otimes X)\\
@VV{\ev_X}V @VV{\omega(\ev_X)}V@VV{\omega(\ev_X)}V\\
\underline{1} @>{\omega}>> \omega(\underline{1}) @>{\id}>>\omega(\underline{1})\\
\end{CD}
\end{equation*}
The evaluation morphism takes $f\otimes x \in \omega(X)^*\otimes\omega(X)$
and evaluates providing $f(x) \in \omega(\underline{1})$ as its output.
Take any $y \in \omega(X)^*$ and $x \in \omega(X)$. Then
\begin{align*}
( y, \mu_X\circ\lambda_X(x)) &= (t_X(y), \lambda_X(x)) =
(\varphi\circ \lambda_{X^*}\circ\varphi^{-1}(y),\lambda_X(x)) =\\
&=\ev_{\omega(X)}\circ(\varphi\otimes \id )\circ(\lambda_{X^*}\otimes \lambda_X)(\varphi^{-1}(y),x)=\\
&=\ev_{\omega(X)}\circ(\varphi\otimes \id )\circ\lambda_{X^*\otimes X}(\varphi^{-1}(y),x)
= (y,x)
\end{align*}
and as a result $\lambda_X$ is injective. Thus, $\lambda_X$ is
also surjective and $\mu_X$ is its inverse. We can show this
differently:
\begin{align*}
( y, \lambda_X\circ\mu_X(x)) =\ev_{\omega(X)}\circ(\varphi\otimes \id )\circ\mu_{X^*\otimes X}(\varphi^{-1}(y),x)
= (y,x)
\end{align*}
as $\mu \in \End^\otimes(\omega)(R).$
\end{proof}

The proof of the following result that we give differs from the similar
one in \cite{Springer}. Our goal
was to provide a correct differential structure.
Also, if one follows our proof in a non-differential case one
finds that it does {\it not} depend on $\Char\K.$ For this
commutative case the change (in comparison to \cite{Springer})
that we make is at the end of the proof where we take the
``generic point''.

\begin{theorem}\label{Recover} We have \begin{enumerate}
\item the algebra $\A$ is a finitely generated $\partial$-$\K$-algebra;
\item there is a surjective $\partial$-$\K$-algebra homomorphism $\Phi : \A \to A$
such that
$$
\Phi\circ a_V = \psi_V
$$
for all differential $G$-modules $V$, where $\psi_V$ is defined in formula~\eqref{SpringersMap};
\item the map $\Phi$ is a $\partial$-$\K$-algebra isomorphism $\A \to A.$
\end{enumerate}
\end{theorem}
\begin{proof}
Let $V \in \Ob(\Rep_G)$. Fix a basis $\{u_1,\ldots,u_n\}$ of $V^*.$
Since $A$ is locally finite and $$\phi_{u_i} : v \mapsto
u_i(r_V(\cdot)\cdot v)$$ is a $G$ module morphism $V \to A$, there
is $W \in \Ob(\Rep_G)$, $W \subset A,$ and $\dim W < \infty$
containing the images of $\phi_{u_i}$ for all $i,$ $1\Le i\Le n.$
According to the proof of Lemma~\ref{Regular} the induced
$G$-morphism $\phi : V \to W^n$ is injective. Hence, the map
$\phi^*$ is surjective and for $u \in V^*$ there exists $t =
(t_1,\ldots,t_n) \in W^n$ such that $u = \phi^*(t)$. We then have
$$
a_V(v\otimes u) = a_V(v\otimes\phi^*((t_1,\ldots,t_n))) =
a_W(\phi(v)\otimes(t_1,\ldots,t_n)).
$$
Thus, the differential algebra $\A$ is generated by the images of the $a_V$ for
$A \supset V \in \Ob(\Rep_G)$ and $\dim V < \infty.$ Let $V$ be such
a $G$-submodule of $A$ which also contains $1$ and a finite set
of generators of $A$ as a $\partial$-$\K$-algebra. The multiplication
on the algebra $A$ defines for any $l \in \mathbb{Z}_{\Ge 1}$ a
surjective $G$-morphism $\phi_l$ from $V^{\otimes l}$ onto some $V(l) \in \Ob(\Rep_G)$
with $\dim V(l) < \infty.$ We then have $V(l) \subset V(l+1),$ because
$1 \in V$, and $A = \bigcup_{l \Ge 1}V(l).$

Consider any $\Ob(\Rep_G) \ni W \subset A$ with $\dim W < \infty.$ There
exists $l \in \mathbb{Z}_{\Ge 1}$ such that $W \subset V(l).$ Since
$\phi_l$ is surjective, we have
$$
\Image a_W\subset \Image a_{V(l)}\subset \Image a_{V^{\otimes l}}.
$$
Because of the multiplication structure of $\A$ the set $\Image a_{V^{\otimes l}}$
lies in the $\partial$-$\K$-subalgebra generated by $\Image a_V$. Hence, this subalgebra is the
whole $\A.$

The homomorphism $\Phi$ of the second statement is constructed as
follows. We take an element $a_V(v\otimes u)$ and map it to $\psi_V(v,u)$
for all $V \in \Ob(\Rep_G).$ Since
\begin{align*}
\Phi\left(a_V(v\otimes u)\cdot
a_W(w\otimes t)\right) &= \Phi\left(a_{V\otimes
W}((v\otimes w)\otimes(u\otimes t))\right)=\\
&=\psi_{V\otimes W}(v\otimes w, u\otimes t) =\\
&= m\circ\left(\psi_V\otimes \psi_W\right)\circ c(v\otimes w,u\otimes t) =\\
&= m\circ\left(\psi_V\otimes
\psi_W\right)(v\otimes u,w\otimes t) = \\
&=m\left(\Phi(a_V(v\otimes u)),
\Phi(a_W(w\otimes t))\right),
\end{align*}
the map $\Phi$ is a $\K$-algebra homomorphism. Let us show that
it is differential. From Lemma~\ref{PhiLemma} we have:
\begin{align*}
\Phi\left(\partial \left(a_V
(v\otimes u\right)\right)
&= \Phi\left(a_{V^{(1)}}
(\partial v\otimes F(u)\right) = \psi_{V^{(1)}}(\partial v\otimes F(u)) =\\
&=\partial\left(\psi_V(v\otimes u)\right)
= \partial\left(\Phi\left(a_V
(v\otimes u)\right)\right).
\end{align*}

We now show the last statement.
Let $B$ be a $\partial$-$\K$-algebra.
Consider a point $\xi \in \Hom_{\K[\partial ]}(\A,B)$ and $V \in \Ob(\Rep_G).$ Fix bases $\{v_i\}$ and $\{u_j\}$ of $V$ and
$V^*$, respectively. There is an endomorphism $\lambda_V$ of $V\otimes B$
such that
$$
\langle\lambda_V(v_i),u_j\rangle = u_j(\lambda_V(v_i)) = \xi(a_V(v_i\otimes u_j)).
$$
We show now
that $(\lambda_V\:|\: V \in \Ob(\Rep_G))$ satisfies the conditions
of Theorem~\ref{TannakaTheorem}.

Let $V, W\in \Ob(\Rep_G)$. Then we have
\begin{align*}
\langle\lambda_{V\otimes W}(v_i\otimes w_j), u_r\otimes t_l\rangle &=
\xi(a_{V\otimes W}((v_i\otimes w_j)\otimes(u_r\otimes t_l))) =\\
&=\xi(a_V(v_i\otimes u_r)\cdot a_W(w_j\otimes t_l)) =\\
&=\xi(a_V(v_i\otimes u_r))\cdot \xi(a_W(w_j\otimes t_l))=\\
&=\langle\lambda_V(v_i),u_r\rangle\cdot\langle\lambda_W(w_j),t_l\rangle =\\
&=\langle(\lambda_V\otimes\lambda_W)(v_i\otimes w_j),u_r\otimes t_l\rangle.
\end{align*}
Hence, $\lambda_{V\otimes W} = \lambda_V\otimes\lambda_W.$ Since
$a_{\1}$ is the identity in $\A$ and $\xi(1_{\A}) = 1$, we have
$\lambda_{\1}$ is the identity. Let us show the functoriality of
$(\lambda_V).$ For a $G$-equivariant map $\phi : V\to W$ we have
\begin{align*}
\langle(\lambda_W\circ \phi)(v_i),t_j\rangle &=\xi(a_W(\phi(v_i)\otimes t_j)) = \xi(a_V(v_i\otimes \phi^*(t_j)))=\\
&= \langle\lambda_V(v_i),\phi^*(t_j)\rangle =
\langle(\phi\circ\lambda_V)(v_i),t_j\rangle.
\end{align*}
Hence, $\lambda_W\circ \phi = \phi\circ\lambda_V.$ Finally, since
$\{v_i\}$ and $\{u_j\}$ are dual to each other, we have:
\begin{align*}
\langle\partial\circ\lambda_V(v_i),u_j\rangle &= \partial(u_j(\lambda_V(v_i)))=
\partial\circ\xi(a_V(v_i\otimes u_j))=\\
 &=\xi(\partial a_V(v_i\otimes u_j))=\xi(a_{V^{(1)}}((\partial v_i) \otimes F(u_j))) = \langle\lambda_{V^{(1)}}(\partial v_i),F(u_j)\rangle.
\end{align*}
Moreover, let $\lambda_V(v_i) = \sum c_{ik}v_k$. Due
to~\eqref{NewDiffEq} we have:
\begin{align*}
\langle\partial\circ\lambda_V(v_i),(\partial v_j)^*\rangle &= (\partial v_j)^*(\partial(c_{ik}v_k)) = (\partial v_j)^*(\partial(c_{ik})v_k+c_{ik}\partial v_k) =\\
&= c_{ik} =\xi(a_V(v_i\otimes v_j^*))= \xi\left(a_{V^{(1)}}\left(\partial v_i\otimes \varphi^*\left(v_j^*\right)\right)\right)=\\
&= \xi(a_{V^{(1)}}(\partial v_i\otimes (\partial v_j)^*))=\langle
\lambda_{V^{(1)}}(\partial v_i),(\partial v_j)^*\rangle.
\end{align*}
We then conclude that $\lambda \in \End^{\otimes,\partial}(\omega).$ Since the category $\Rep_G$
is rigid, $\lambda \in \AutOP(\omega)$ by Lemma~\ref{RigidLemma}.

Take $B = \A$ and the generic point $\xi = \id_{\A}.$ By Theorem~\ref{TannakaTheorem}
there exists $x \in \Hom_{\K[\partial ]}(A,\A) = G(\A)$ such that
$\lambda_V = r_V(x)$ for all $V \in \Ob(\Rep_G).$ We have
\begin{align*}
x\circ\Phi\circ a_V(v_i\otimes u_j)) &= x\circ\psi_V(v_i\otimes u_j) = \langle r_V(x)v_i,u_j\rangle=\\
&=\langle\lambda_V(v_i),u_j\rangle = \xi(a_V(v_i\otimes u_j)).
\end{align*}
Hence,
$x\circ\Phi = \xi = \id_{\A}.$ This implies that $\Phi$ is injective. Since
$\Phi$ is also surjective, we obtain that $\Phi : \A \to A$ is a
$\partial$-$\K$-algebra isomorphism.
\end{proof}

\subsection{Recovering $\Delta$ and $S$}\label{RecoverDeltaS}
We provide a differential Hopf algebra structure to $\A.$
Let $V \in \Ob(\Rep_G)$ and $\{v_i\}$ be its basis with the dual basis
$\{u_j\}$ of $V^*.$
Recall the $\K$-linear map
$$
\tilde{\Delta} : \A \to \A\otimes\A,\quad a_V(v\otimes u)
\mapsto\sum a_V(v_i\otimes u)\otimes a_V(v\otimes u_i).
$$
\begin{lemma} The map $\tilde{\Delta}$ is a $\partial$-$\K$-algebra
homomorphism and is a comultiplication.
\end{lemma}
\begin{proof}
We first check the basis independence. Let $\{e_1,\ldots,e_n\}$ be
another basis for $V$ and $\{f_1,\ldots,f_n\}$ be its dual. Hence,
there exits a matrix $C = (c_{ij}) \in \GL_n(\K)$ such that $v_i =
\sum e_j c_{ji}.$ We then have:
\begin{align*}
\sum_{i=1}^n a_V(v_i\otimes u)\otimes a_V(v\otimes u_i) &= \sum_{i=1}^n a_V\left(\sum_{j=1}^n e_jc_{ji}\otimes u\right)\otimes a_V(v\otimes u_i)=\\
&= \sum_{j=1}^n a_V\left(e_j\otimes u\right)\otimes a_V\left(v\otimes\sum_{i=1}^n c_{ji} u_i\right)=\\
&=\sum_{j=1}^n a_V\left(e_j\otimes u\right)\otimes a_V(v\otimes f_j).
\end{align*}
We check that it is an algebra homomorphism. We have:
\begin{align*}
&\tilde{\Delta}(a_V(v\otimes u)\cdot a_W(w\otimes t)) =
 \tilde{\Delta}(a_{V\otimes W}((v\otimes w)\otimes(u\otimes t)))=\\
&= \sum_{i,j} a_{V\otimes W}((v_i\otimes w_j)\otimes(u\otimes t))\otimes a_{V\otimes W}((v\otimes w)\otimes(u_i\otimes t_j))=\\
&=\sum_{i,j} (a_V(v_i\otimes u)\otimes a_V(v\otimes u_i))
\cdot (a_W(w_j\otimes t)\otimes a_W(w\otimes t_j))=\\
&=\left(\sum_i a_V(v_i\otimes u)\otimes a_V(v\otimes u_i)\right)
\cdot\left(\sum_j a_W(w_j\otimes t)\otimes a_W(w\otimes t_j)\right)=\\
&=\tilde{\Delta}(a_V(v\otimes u))\cdot\tilde{\Delta}(a_W(w\otimes t)).
\end{align*}
Moreover, due to identity~\eqref{NewDiffEq} and the imbedding
$$V \mapsto V^{(1)},\quad v \mapsto 1\otimes v,\ v \in V,$$
we have
\begin{align*}
\partial\circ\tilde{\Delta}(a_V(v\otimes u)) = &\sum \left(a_{V^{(1)}}(\partial v_i\otimes F(u))\otimes a_V(v\otimes u_i)+\right.\\ &\left.+a_V(v_i\otimes u)\otimes a_{V^{(1)}}(\partial v\otimes F(u_i))\right) =\\
=&\sum \left(a_{V^{(1)}}(\partial v_i\otimes F(u))\otimes a_{V^{(1)}}(\partial v\otimes(\partial v_i)^*)+\right.\\
 &\left.+a_{V^{(1)}}(v_i\otimes u)\otimes a_{V^{(1)}}(\partial v\otimes F(u_i))\right) =\\
= &\tilde{\Delta}(\partial(a_V(v\otimes u))).
\end{align*}

We finally show the coassociativity:
\begin{align*}
\left(\tilde{\Delta}\otimes \id\right)\circ\tilde{\Delta}(a_V(v\otimes u)) &= \sum_i \tilde{\Delta}(a_V(v_i\otimes u))\otimes a_V(v\otimes u_i) =\\
&=\sum_i \left(\sum_j a_V(v_j \otimes u)\otimes a_V(v_i\otimes u_j)\right)\otimes a_V(v\otimes u_i)=\\
&=\sum_i a_V(v_i \otimes u)\otimes\left(\sum_j a_V(v_j\otimes u_i)\otimes a_V(v\otimes u_j)\right)=\\
&=\sum_i a_V(v_i \otimes u)\otimes\tilde{\Delta}(a_V(v\otimes u_i))=\\
&=\left(\id\otimes\tilde{\Delta}\right)\circ\tilde{\Delta}(a_V(v\otimes u)).
\end{align*}
\end{proof}

\begin{lemma}
The map $\tilde{\varepsilon}: a_V(v\otimes u) \mapsto u(v)$ is a
counit for $\A$ corresponding to the counit of $A.$
\end{lemma}
\begin{proof}
We show that $\tilde{\varepsilon}$ and  $\tilde{\Delta}$ satisfy
$m\circ(\id_{\A}\otimes\tilde{\varepsilon})\circ\tilde{\Delta} =
\id_{\A}.$ We have:
\begin{align*}
m\circ(\id_{\A}\otimes\tilde{\varepsilon})\circ\tilde{\Delta}(a_V(v\otimes u)) &=
m\circ(\id_{\A}\otimes\tilde{\varepsilon})\left(\sum a_V(v_i\otimes u)\otimes a_V(v\otimes u_i)\right) =\\
&=\sum(a_V(v_i\otimes u)\cdot u_i(v)) = \\
&= a_V\left(\left(\sum u_i(v)\cdot v_i\right)\otimes u\right) = a_V(v\otimes u).
\end{align*}
In addition, $\tilde{\varepsilon}$ is a differential homomorphism.
Indeed,
$$
\tilde{\varepsilon}\left(\partial a_V(v\otimes u)\right) =
\tilde{\varepsilon}\left(a_{V^{(1)}}(\partial v\otimes F(u))\right)
= F(u)(\partial v) = \partial(u(v)) =
\partial(\tilde{\varepsilon}(a_V(v\otimes u))).
$$
Finally, we show that $\Phi$ maps $\tilde{\varepsilon}$ to $\varepsilon.$
\begin{align*}
\varepsilon\circ\Phi(a_V(v\otimes u)) = \Phi(a_V(v\otimes u))(e) = u(r_V(e)\cdot v) = u(v) = \Phi(\tilde{\varepsilon}(a_V(v\otimes u)).
\end{align*}
\end{proof}

\begin{proposition}\label{HopfAlgebraHomomorphism} The map $\Phi : \A \to A$ is a differential Hopf algebra
homomorphism.
\end{proposition}
\begin{proof}
We show that $\tilde{\Delta}$ is mapped to $\Delta.$ We have:
\begin{align*}
&(\Phi\otimes\Phi)\circ\tilde{\Delta}(a_V(v\otimes u))(g_1,g_2) =\\
&= (\Phi\otimes\Phi)\left(\sum_i(a_V(v_i\otimes u)\otimes a_V(v\otimes u_i)\right)(g_1,g_2)=\\
&=\sum_i u(r_V(g_1)\cdot v_i)\cdot u_i(r_V(g_2)\cdot v);\\
&\Delta\circ\Phi(a_V(v\otimes u))(g_1,g_2) = u(r_V(g_1\cdot g_2)\cdot v) =\\
&= u(r_V(g_1)\cdot(r_V(g_2)\cdot v)).
\end{align*}
Let $r_V(g_2)\cdot v = \sum c_jv_j.$ Then
\begin{align*}
\sum_i u(r_V(g_1)\cdot v_i)\cdot u_i(r_V(g_2)\cdot v) &=
\sum_i u(r_V(g_1)\cdot v_i)\cdot u_i\left(\sum_jc_jv_j\right) =\\
&=\sum_i u(r_V(g_1)\cdot v_i)\cdot c_i =\\
&=\sum_i u(r_V(g_1)\cdot c_iv_i) =\\
&=u(r_V(g_1)\cdot (r_V(g_2)\cdot v)).
\end{align*}
\end{proof}

Recall the $\K$-linear map
$$
\tilde{S} : \A \to \A,\quad a_V(v\otimes u)
\mapsto a_{V^*}(u\otimes v).
$$
\begin{lemma} The map $\tilde{S}$ is a $\partial$-$\K$-algebra
homomorphism and together with $\tilde{\Delta}$ gives a
differential Hopf algebra structure on $\A$.
\end{lemma}
\begin{proof}
Let $\{v_1,\ldots, v_n\}$ be a basis of $V$ and $\{v_1^*,\ldots,v_n^*\}$
be its dual.
We show that $\Coinv$ commutes with $\partial:$
\begin{align*}
\partial\left(\tilde{S}(a_V(v_i\otimes v_j^*))\right) &= \partial (a_{V^*}(v_j^*\otimes v_i))=
a_{{(V^*)}^{(1)}}\left(\partial (v_j^*)\otimes F(v_i)\right) =\\
&= a_{{\left(V^{(1)}\right)}^*}(F(v_j^*)\otimes\partial
v_i)=\tilde{S}\left(a_{V^{(1)}}\left(\partial v_i\otimes
F(v_j^*)\right)\right) =\\
&= \tilde{S}(\partial a_V(v_i\otimes v_j^*)),
\end{align*}
where we use the morphism
$$\phi : \left(V^*\right)^{(1)} \to \left(V^{(1)}\right)^*,\quad \partial v_j^* \mapsto F(v_j^*),\ v_j^* \mapsto (\partial v_j)^*$$
commuting with the $G$-action.
Moreover, $\Coinv$ is an algebra homomorphism:
\begin{align*}
\tilde{S}(a_V(v\otimes u)\cdot a_W(w\otimes t)) &= \tilde{S}(a_{V\otimes W}((v\otimes w)\otimes(u\otimes t))) =\\
&= a_{V^*\otimes W^*}((u\otimes t)\otimes(v\otimes w)) =\\
&= a_{V^*}(u\otimes v)\cdot a_{W^*}(t\otimes w)=\\
&= \tilde{S}(a_V(v\otimes u))\cdot\tilde{S}(a_W(w\otimes t)).
\end{align*}
We show that it respects comultiplication:
\begin{align*}
m\circ(\tilde{S}\otimes\id)\circ\tilde{\Delta}(a_V(v\otimes u)) &=
m\circ(\tilde{S}\otimes\id)\left(\sum_ia_V(v_i\otimes u)\otimes a_V(v\otimes u_i)\right) =\\
&=m\left(\sum_ia_{V^*}(u\otimes v_i)\otimes a_V(v\otimes u_i)\right) =\\
&=\sum_ia_{V^*\otimes V}((u\otimes v)\otimes (v_i\otimes u_i)) =\\
&=u(v)\cdot a_I(e\otimes f) = \\
&=\E(a_V(v\otimes u)).
\end{align*}
We have denoted a basis of the trivial representation $\1$ by $\{e\}$
and its dual by $\{f\}$.
We have also used the morphism $V^*\otimes V \to \1$
mapping $u\otimes v$ to $u(v)$.
\end{proof}

\section{Partial differential case}\label{ManyParameters}
We have only used elementary ring theoretic properties  of $\K[\partial]$
and none of its special properties as a left and right Euclidean domain.
In particular, all statements concerning recovering the differential Hopf algebra from
representations hold true in the partial case.
We just restate Definition~\ref{SeqCategory} for the case of
several commuting differentiations.
\begin{definition}\label{PartialSeqCategory} The category $\Seq_{\K}({\partial_1,\ldots,\partial_m})$ over a $\{\partial_1,\ldots,\partial_m\}$-field $\K$ is
the category of finite dimensional vector spaces together
with the usual operations $\otimes,$ $\oplus,$ $^*$, and additional
differentiation functors $$
\partial_1^{p_1}\cdot\ldots\cdot\partial_m^{p_m} : V \mapsto \K[\partial_1,\ldots,\partial_m]_{\Le (p_1,\ldots,p_m)}\otimes V
$$
for all $m$-tuples $(p_1,\ldots,p_m) \in \left(\mathbb{Z}_{\Ge 0}\right)^m.$
\end{definition}

\section{Conclusions}
The results of the previous section allow us to recover a differential algebraic group from the category
of its finite dimensional differential representations. From Proposition~\ref{GenerateRepresentations}
this category can be generated by one faithful representation of the group applying
certain operations of linear algebra and the prolongation functor we introduced in this paper.

\section{Acknowledgements}
The author is highly grateful to his advisor Michael Singer,
to Bojko Bakalov, Pierre Deligne, and Daniel Bertrand
for extremely helpful comments and support. Also, the author
thanks the participants of Kolchin's Seminar in New York for
their important suggestions. The author appreciates the 
detailed comments of the referees very much.

\end{document}